\documentclass{amsart}
\usepackage{amsopn}
\usepackage{amssymb, amscd}

\newcommand{\nc}{\newcommand}

\nc{\vg}{\mathfrak{v} } \nc{\wg}{\mathfrak{w} } \nc{\zg}{\mathfrak{z} } \nc{\ngo}{\mathfrak{n} }
\nc{\kg}{\mathfrak{k} } \nc{\mg}{\mathfrak{m} } \nc{\bg}{\mathfrak{b} } \nc{\ggo}{\mathfrak{g} }
\nc{\ggob}{\overline{\mathfrak{g}} } \nc{\sog}{\mathfrak{so} } \nc{\sug}{\mathfrak{su} } \nc{\spg}{\mathfrak{sp}
} \nc{\slg}{\mathfrak{sl} } \nc{\glg}{\mathfrak{gl} } \nc{\cg}{\mathfrak{c} } \nc{\rg}{\mathfrak{r} }
\nc{\hg}{\mathfrak{h} } \nc{\tg}{\mathfrak{t} } \nc{\ug}{\mathfrak{u} } \nc{\dg}{\mathfrak{d} }
\nc{\ag}{\mathfrak{a} } \nc{\pg}{\mathfrak{p} } \nc{\sg}{\mathfrak{s} } \nc{\pca}{\mathcal{P}}
\nc{\nca}{\mathcal{N}} \nc{\lca}{\mathcal{L}} \nc{\bca}{\mathcal{B}} \nc{\oca}{\mathcal{O}}
\nc{\mca}{\mathcal{M}} \nc{\tca}{\mathcal{T}} \nc{\aca}{\mathcal{A}} \nc{\vp}{\varphi} \nc{\ddt}{\frac{{\rm
d}}{{\rm d}t}} \nc{\im}{\mathtt{i}}
\renewcommand{\Re}{{\rm Re}}
\renewcommand{\Im}{{\rm Im}}

\nc{\SO}{{\mathrm SO}} \nc{\Spe}{{\mathrm Sp}} \nc{\Sl}{{\mathrm SL}} \nc{\SU}{{\mathrm SU}} \nc{\Or}{{\mathrm
O}} \nc{\U}{{\mathrm U}} \nc{\Gl}{{\mathrm GL}} \nc{\Se}{{\mathrm S}} \nc{\Cl}{{\mathrm Cl}} \nc{\Spin}{{\mathrm
Spin}} \nc{\Pin}{{\mathrm Pin}}

\nc{\RR}{{\Bbb R}} \nc{\HH}{{\Bbb H}} \nc{\CC}{{\Bbb C}} \nc{\ZZ}{{\Bbb Z}} \nc{\FF}{{\Bbb F}} \nc{\NN}{{\Bbb
N}} \nc{\QQ}{{\Bbb Q}} \nc{\PP}{{\Bbb P}}

\nc{\vs}{\vspace{.2cm}} \nc{\vsp}{\vspace{1cm}} \nc{\ip}{\langle\cdot,\cdot\rangle} \nc{\la}{\langle}
\nc{\ra}{\rangle} \nc{\unm}{\frac{1}{2}} \nc{\unc}{\frac{1}{4}} \nc{\unt}{\frac{1}{3}} \nc{\und}{\frac{1}{16}}
\nc{\no}{\vs\noindent} \nc{\lam}{\Lambda^2\ggo^*\otimes\ggo} \nc{\Fres}{F|_{S}} \nc{\tangz}{{\rm T}^{\rm Zar}}
\nc{\nor}{{\sf n}} \nc{\eigen}{(k_1<...<k_r;d_1,...,d_r)} \nc{\eigencero}{(0<k_2<...<k_r;d_1,...,d_r)}
\nc{\mum}{/\!\!/} \nc{\kir}{/\!\!/\!\!/}

\nc{\He}{\operatorname{Hess}} \nc{\ad}{\operatorname{ad}} \nc{\Ad}{\operatorname{Ad}}
\nc{\rank}{\operatorname{rank}} \nc{\Irr}{\operatorname{Irr}} \nc{\End}{\operatorname{End}}
\nc{\Aut}{\operatorname{Aut}} \nc{\Inn}{\operatorname{Inn}} \nc{\Der}{\operatorname{Der}}
\nc{\Ker}{\operatorname{Ker}} \nc{\Iso}{\operatorname{I}} \nc{\Diff}{\operatorname{D}}
\nc{\Lie}{\operatorname{L}} \nc{\tr}{\operatorname{tr}} \nc{\dif}{\operatorname{d}}
\nc{\sen}{\operatorname{sen}} \nc{\modu}{\operatorname{mod}} \nc{\Ric}{\operatorname{Ric}}
\nc{\sym}{\operatorname{sym}} \nc{\sca}{\operatorname{sc}} \nc{\scalar}{{\sf s}} \nc{\grad}{\operatorname{grad}}
\nc{\ricci}{\operatorname{ric}} \nc{\Rin}{\operatorname{M}} \nc{\Le}{\operatorname{L}}
\nc{\tang}{\operatorname{T}} \nc{\level}{\operatorname{level}} \nc{\rad}{\operatorname{r}}
\nc{\abel}{\operatorname{ab}} \nc{\Spec}{\operatorname{Spec}}

\newtheorem{theorem}{Theorem}[section]
\newtheorem{proposition}[theorem]{Proposition}
\newtheorem{corollary}[theorem]{Corollary}
\newtheorem{lemma}[theorem]{Lemma}
\newtheorem{definition}[theorem]{Definition}
\newtheorem{remark}[theorem]{Remark}

\newtheorem{example}[theorem]{Example}

\title{On the moment map for the variety of Lie algebras}

\author{Jorge Lauret}

\address{Department of Mathematics, Yale University,
10 Hillhouse Box 208283 New Haven, CT 06520 USA (current affiliation: postdoctoral fellowship)}
\email{jorge.lauret@yale.edu}

\address{FaMAF and CIEM, Universidad Nacional de C\'ordoba, 5000 C\'ordoba, Argentina}
\email{lauret@mate.uncor.edu}

\thanks{
Supported by CONICET and Guggenheim Foundation fellowships, and grants from FONCyT and SeCyT UNC (Argentina).}

\begin{document}

\maketitle

\begin{abstract}
We consider the moment map $m:\PP V_n\mapsto\im\ug(n)$ for the action of $\Gl(n)$ on
$V_n=\Lambda^2(\CC^n)^*\otimes\CC^n$.  The critical points of the functional $F_n=||m||^2:\PP V_n\mapsto\RR$ are
studied, in order to understand the stratification of $L_n\subset\PP V_n$ defined by the negative gradient flow
of $F_n$, where $L_n$ is the algebraic subset of all Lie algebras.  We obtain a description of the critical
points which lie in $L_n$ in terms of those which are nilpotent, as well as the minima and maxima of
$F_n:L_n\mapsto\RR$.  A characterization of the critical points modulo isomorphism, as the union of categorical
quotients of suitable actions is considered, and some applications to the study of $L_n$ are given.
\end{abstract}

\section{Introduction}\label{intro}

The space of all complex Lie algebras of a given dimension $n$ can be naturally identified with the set $\lca_n$
of all Lie brackets on $\CC^n$.  Since the Jacobi identity is determined by polynomial conditions, $\lca_n$ is
an algebraic subset of the vector space $V_n$ of skew-symmetric bilinear maps from $\CC^n\times\CC^n$ to
$\CC^n$. The isomorphism class of a Lie algebra $\mu\in\lca_n$ is then given by the orbit $\Gl(n).\mu$ under the
`change of basis' action of $\Gl(n)$ on $V_n$.  This action is very unpleasant from the point of view of
invariant theory since any $\mu\in V_n$ is unstable (i.e. $0\in\overline{\Gl(n).\mu}$), which makes very
difficult the study of the quotient space $\lca_n/\Gl(n)$ parameterizing Lie algebras up to isomorphism.

Nevertheless, F. Kirwan \cite{Krw1} and L. Ness \cite{Ns} have showed that the moment map for an action can be
used to study the orbit space of the null-cone (set of unstable vectors).  Let $m:\PP V\mapsto\im\kg$ be the
moment map for the action of a complex reductive Lie group $G$ with maximal compact subgroup $K$ on a vector
space $V$.  An orbit $G.v$ is closed if and only if $G.[v]$ meets $m^{-1}(0)$ and the intersection is a single
$K$-orbit (see \cite{KmNs}), where $[v]$ denotes the class of $v$ in the projective space $\PP V$.  Let
$X\subset \PP V$ be a $G$-invariant projective algebraic variety.  The so called categorical quotient $X\mum G$
parameterizing closed orbits is homeomorphic to $(X\cap m^{-1}(0))/K$, which is precisely the symplectic
reduction when $X$ is nonsingular.

It is proved in \cite{Krw1} and \cite{Ns} that the remaining critical points of $F=||m||^2:\PP V\mapsto\RR$
(i.e. such that $F(x)>0$) are all in the null-cone and in some sense, their orbits play the same role that the
closed orbits in the set of semistable points (i.e. $0\notin\overline{G.v}$).  For instance, if $C$ denotes the
set of critical points of $F$ then $C\cap G.[v]$ is either empty or a single $K$-orbit, and they are minima of
$F|_{G.[v]}$.  It is then natural to consider a wider quotient $X\kir G$, which shall be called {\it Kirwan-Ness
quotient}, given by
$$
X\kir G=(X\cap C)/K.
$$
Recall that if the action contains the homotheties $\{ v\mapsto tv:t\in\CC^*\}$, then we may also define the
Kirwan-Ness quotient for the action of $G$ on an algebraic $G$-variety $X\subset V$ by
$$
X\kir G=\{ v\in X:[v]\in C\}/\CC^*K,
$$
and clearly $X\kir G=\pi(X)\kir G$, where $\pi:V\setminus\{ 0\}\mapsto\PP V$ is the usual projection map.  This
new quotient $X\kir G$ is not a projective algebraic variety (not even Hausdorff) as in the case of $X\mum G$,
but nevertheless, its topology is not so wild.  Indeed, $X\kir G$ can be decomposed as a disjoint union of
projective algebraic varieties with a respectable frontier property, coming from the stratification defined by
the negative gradient flow of $F$ (see Section \ref{ness}).

Let us now go back to the variety of Lie algebras.  We note that $\lca_n\mum\Gl(n)$ consists of only one point,
the abelian Lie algebra, and we will show that to consider $\Sl(n)$-orbits does not help much, indeed
$\lca_n\mum\Sl(n)=$ semisimple Lie algebras (plus the abelian Lie algebra).  The aim of this paper is to
initiate the study of the Kirwan-Ness quotient $\lca_n\kir\Gl(n)$, or $L_n\kir\Gl(n)$, where $L_n=\pi(\lca_n)$.
How special are the Lie algebras which are isomorphic to a critical point of $F_n:\PP V_n\mapsto\RR$?.  The main
result is that the study of $\lca_n\kir\Gl(n)$ reduces essentially to the understanding of $\nca_n\kir\Gl(n)$,
where $\nca_n$ is the subvariety of nilpotent Lie algebras.  Our real goal is however the possible applications
of this `moment map' approach to the study of degenerations, rigidity and of the irreducible components of
$\lca_n$.

The moment map is in our case a $\U(n)$-equivariant function $m:\PP V_n\mapsto\im\ug(n)$, where $\im\ug(n)$
denotes the space of hermitian matrices.  It comes from the Hamiltonian action of $\U(n)$ on the symplectic
manifold $\PP V_n$.  One fixes an inner product $\ip$ on $\CC^n$ and considers the corresponding
$\U(n)$-invariant inner products on $V_n$ and $\im\ug(n)$, respectively, naturally associated with $\ip$.  For
each $\mu\in V_n$ with $||\mu||=1$, $m([\mu])$ is defined as the derivative at the identity of the function
$\Gl(n)\mapsto\RR$, $g\mapsto ||g.\mu||^2$.  It will be proved that
$$
m([\mu])= \Rin_{\mu}:=-4\displaystyle{\sum_{i}}(\ad_{\mu}{X_i})^*\ad_{\mu}{X_i}
+2\displaystyle{\sum_{i}}\ad_{\mu}{X_i}(\ad_{\mu}{X_i})^*, \qquad ||\mu||=1,
$$
where $\{ X_1,...,X_n\}$ is any orthonormal basis of $\CC^n$.  We now consider the functional $F_n:\PP
V_n\mapsto\RR$ given by $F_n([\mu])=||m([\mu])||^2=\tr{\Rin_{\mu}^2}$.  Some of the remarkable properties of
$F_n$, its gradient flow and their critical points can be summarized as follows (see Section \ref{bil}):

\begin{itemize}
\item $\grad(F_n)_{[\mu]}=-8\pi_{\ast}\delta_{\mu}(\Rin_{\mu})$,  $||\mu||=1$, where $\delta_\mu:\glg(n)\mapsto
V_n$ coincides with the coboundary operator relative to adjoint cohomology when $\mu$ is a Lie algebra, and
$\pi_{\ast}$ denotes the derivative of $\pi:V_n\setminus\{ 0\}\mapsto \PP V_n$.

\item The gradient of $F_n$ is then always tangent to the $\Gl(n)$-orbits, and so $[\mu]$ is a critical point of
$F_n$ if and only if it is such for $F_n|_{\Gl(n).[\mu]}$.  Thus, the negative gradient flow of $F_n$ stays in
the orbit of the starting point $[\lambda]$ and therefore it gives rise a distinguished degeneration from
$\lambda$ to a critical point $\mu$ (i.e. $\mu\in\overline{\Gl(n).\lambda}$).

\item If $[\mu]$ is a critical point of $F_n$ then $F_n|_{\Gl(n).[\mu]}$ attains its minimum value at $[\mu]$,
and any other critical point in $\Gl(n).[\mu]$ belongs to $\U(n).[\mu]$.

\item $[\mu]\in\PP V_n$ is a critical point of $F_n$ if and only if $M_{\mu}=cI+D$ for some $c\in\RR$ and
$D\in\Der(\mu)$. In this case, a positive multiple of $D$ has nonnegative integer eigenvalues $k_i$, say with
multiplicities $d_i$, and $c$ can be computed in terms of these integers.  The critical point $[\mu]$ is said to
be of {\it type} $\alpha=(k_1<...<k_r;d_1,...,d_r)$.

\item Any critical point $[\mu]$ admits a $\ZZ_{\geq 0}$-gradation.  If $\mu\in\lca_n$ then $k_1>0$ if and only
if $\mu$ is nilpotent, and so any nilpotent Lie algebra which is a critical point of $F_n$ is $\NN$-graded.

\item There are finitely many types of critical points, say $\alpha_1,...,\alpha_s$.

\item If $C_{\alpha}$ is the set of critical points of $F_n$ of type $\alpha$, then the quotient space
$C_{\alpha}/\U(n)=\Gl(n).C_{\alpha}/\Gl(n)$ is homeomorphic to the categorical quotient of a suitable action,
and so it is a projective algebraic variety (see \cite{MFK}).

\item Let $S_{\alpha}\subset\PP V_n$ be the set of all the points which are carried by the negative gradient
flow of $F_n$ to a critical point of type $\alpha$.  Then
$$
\PP V_n=S_{\alpha_1}\cup...\cup S_{\alpha_s}
$$
determines a stratification of $\PP V_n$, for which each stratum $S_{\alpha_i}$ is locally closed, irreducible
and nonsingular.
\end{itemize}

We are interested in the stratification of $L_n$ given by
$$
L_n=(S_{\alpha_1}\cap L_n)\cup...\cup (S_{\alpha_s}\cap L_n),
$$
and consequently in the critical points of $F_n$ which lie in $L_n$ (see Section \ref{critvar}).  The
Kirwan-Ness quotient admits a decomposition
$$
L_n\kir\Gl(n)=X_{\alpha_1}\cup...\cup X_{\alpha_s} \qquad \mbox{(disjoint union)},
$$
where each $X_{\alpha_i}$ is homeomorphic to $C_{\alpha_i}/\U(n)$ and the following frontier property holds:
there is a partial order on the indexing set $\{\alpha_1,...,\alpha_s\}$ such that
$$
\overline{X}_{\alpha}\subset X_{\alpha}\cup\bigcup_{\beta>\alpha}X_{\beta}.
$$
If $\alpha=(0;n)$ then $X_{\alpha}$ consists of finitely many points: the semisimple Lie algebras of dimension
$n$.

We first study extremal points of $F_n:L_n\mapsto\RR$, proving that the minimum value is attained at semisimple
Lie algebras and the maximum value at the direct sum of the $3$-dimensional Heisenberg Lie algebra and the
abelian algebra.

As expected, several strong compatibility properties between a Lie bracket $\mu$ and the fixed inner product
$\ip$ are necessary in order for $[\mu]$ to be a critical point of $F_n$.  For instance, there must exist an
orthonormal decomposition $\CC^n=\hg\oplus\ag\oplus\ngo$, with $\hg$ a semisimple Lie subalgebra of $\mu$, $\ag$
abelian, $\mu(\hg,\ag)=0$ and $\ngo$ the nilradical of $\mu$, such that the adjoint action of the reductive part
$\hg\oplus\ag$ of $\mu$ on the underlying inner product space $(\CC^n,\ip)$ is as nice as possible.  More
precisely,
\begin{itemize}
\item[(i)] $\ad_{\mu}{A}$ is a normal operator (and so semisimple) for every $A\in\ag$.

\item[(ii)] The real subalgebra $\kg=\{ A\in\hg:(\ad_{\mu}A)^*=-\ad_{\mu}A\}$ is a maximal compact subalgebra of
$\hg$, that is, $\hg=\kg+\im\kg$.

\item[(iii)] The hermitian inner product $\ip$ on $\hg\oplus\ag$ is given by
$$
\la A,B\ra=-\frac{4}{c_{\mu}}\left(\unm\tr{\ad_{\mu}{A}
(\ad_{\mu}{B})^*}|_{\hg}+\tr{\ad_{\mu}{A}(\ad_{\mu}{B})^*}|_{\ngo}\right), \qquad A,B\in\hg\oplus\ag.
$$
\end{itemize}
We have obtained the following characterization of Lie algebras which are critical points of $F_n$ (see Theorem
\ref{car} for a precise statement).

\begin{theorem}
$[\mu]\in L_n$ is a critical point of $F_n$ if and only if there exists an orthonormal decomposition
$\CC^n=\rg\oplus\ngo$, where $\rg$ is a reductive Lie subalgebra of $\mu$ with the properties that
$(\ad_{\mu}{A})^*\in\Der(\mu)$ for any $A\in\rg$, $\ngo$ is the nilradical of $\mu$ and $\mu|_{\ngo\times\ngo}$
is also a critical point of the corresponding $F_m$, $m=\dim{\ngo}$.
\end{theorem}

Thus, the classification of critical points of $F_n$ in $L_n$ reduces to the determination of those which are
nilpotent.  There is an intriguing interplay between nilpotent critical points and Riemannian geometry (see
Remark \ref{riem}).

In Section \ref{abelian}, the closed subset $\aca\subset\lca_n$ of Lie algebras having a codimension one abelian
ideal is considered, in order to exemplify most of the notions studied in this paper.  Finally, in Section
\ref{critclosed}, we consider for each type $\alpha=\eigen$ the action of the reductive Lie subgroup
$$
\tilde{G}_{\alpha}=\left\{ g\in
\Gl(d_1)\times...\times\Gl(d_r):\prod_{i=1}^{r}(\det{g_i})^{k_i}=\det{g}=1\right\}
$$
on
$$
V_{\alpha}=\{\mu\in V_n:D_{\alpha}\in\Der(\mu)\},
$$
where $D_{\alpha}$ is the diagonal matrix with entries $k_i$ and multiplicities $d_i$.  If $\tilde{m}:\PP
V_{\alpha}\mapsto\im\tilde{\kg}_{\alpha}$ denotes the moment map for this action then $\tilde{m}^{-1}(0)=\PP
V_{\alpha}\cap C_{\alpha}$ and the semistable points $\PP V_{\alpha}^{\rm ss}=\PP V_{\alpha}\cap S_{\alpha}$.
Thus the categorical quotient $\PP V_{\alpha}//\tilde{G}_{\alpha}$ coincides with $\PP V_{\alpha}\cap
C_{\alpha}/\tilde{K}_{\alpha}= C_{\alpha}/\U(n)$, which parameterizes the critical points of type $\alpha$,
modulo isomorphism.

As an application, we study the $\Sl(n)$-action on $\lca_n$, obtaining that the semistable points are precisely
the semisimple Lie algebras and moreover, we show that $\Sl(n).\mu$ is closed if and only if $\mu$ is
semisimple.  A new proof of the rigidity (i.e. $\Gl(n).\mu$ open in $\lca_n$) of semisimple Lie algebras is also
obtained.  We finally give a complete description of the case $n=4$.

\no {\it Acknowledgements.}  I wish to thank Fritz Grunewald for very helpful comments on a first version of
this paper.

\section{A stratification of the null-cone}\label{ness}

Let $G$ be a complex reductive Lie group acting on a finite dimensional complex vector space $V$, and let
$X\subset V$ be a $G$-variety, that is, an algebraic variety which is $G$-invariant.  The main problem of
geometric invariant theory is to understand the orbit space of the action of $G$ on $X$, parameterized by the
quotient $X/G$ (we refer to \cite{PV} for further information).  The standard quotient topology of $X/G$ can be
very ugly, for instance, if $x$ degenerates to $y$ (i.e. $y\in\overline{G.x}$) and $G.x\ne G.y$ then they can
not be separated by $G$-invariant open neighborhoods and so $X/G$ is usually non-Hausdorff.

In order to avoid this problem one may consider a smaller quotient $X\mum G$ parame\-tri\-zing only closed
orbits.  D. Mumford \cite{MFK} proved that $X\mum G$ is again an algebraic variety; indeed, $X\mum
G=\Spec(\CC[X]^G)$, that is, $X\mum G$ is the algebraic variety with coordinate ring $\CC[X]^G$, the
$G$-invariant polynomials on $X$.  Consider $q:X\mapsto X\mum G$ the morphism of algebraic varieties determined
by the inclusion $\CC[X]^G\mapsto\CC[X]$.  $X\mum G$ is called the {\it categorical quotient} for the action of
$G$ on $X$ since it satisfies the following universal property in the category of all algebraic varieties: for
any morphism of algebraic varieties $\alpha:X\mapsto Y$ that is constant on $G$-orbits there exists a unique
morphism $\beta:X\mum G\mapsto Y$ such that $\alpha=\beta\circ q$.  The uniqueness of an object with such a
property is clear from the definition.  D. Luna proved that actually, $X\mum G$ satisfies the same universal
property in the category of all Hausdorff topological spaces.   Recall that the usual quotient $X/G$ would be
the categorical quotient in the category of topological spaces.

The price to pay for a Hausdorff quotient is that in some cases $X\mum G$ classifies only a very few orbits.  If
we have for example that the homotheties $\{v\mapsto tv:t\in\CC^*\}$ are contained in the action of $G$ on $V$
then $\{ 0\}$ will be the unique closed orbit, and hence $X\mum G$ will consist of only one point: too
expensive.  This is also clear from the fact that in such a case $\CC[V]^G=$ constant polynomials.

The aim of this section is to describe a wider quotient $X\kir G$, which of course is non-Hausdorff but it still
satisfies some properties very similar in spirit to those of $X\mum G$.  We call $X\kir G$ the {\it Kirwan-Ness
quotient} for the action of $G$ on $X$ since its definition comes from some remarkable properties of the moment
map for an action proved independently by L. Ness \cite{Ns} and F. Kirwan \cite{Krw1}, which we now overview.

Assume $V$ is endowed with an hermitian inner product $\ip$ which is invariant under the action of the maximal
compact subgroup $K\subset G$.  For each $v\in V$ define
$$
\rho_v:G\mapsto\RR, \qquad \rho_v(g)=||g.v||^2=\la g.v,g.v\ra.
$$
Let $(\dif\rho_v)_e:\ggo\mapsto\RR$ denote the differential of $\rho_v$ at the identity $e$ of $G$.  It follows
from the $K$-invariance of $\ip$ that $(\dif\rho_v)_e$ vanishes on $\kg$, and so we may view
$(\dif\rho_v)_e\in\ggo^*/\kg^*$, where $\ggo^*$ and $\kg^*$ are the vector spaces of real-valued functionals on
the Lie algebras of $G$ and $K$ respectively.  Since $G$ is reductive $\ggo=\kg+\im\kg$, thus we may define a
function
\begin{equation}\label{moment}
m:\PP V\mapsto\im\kg, \qquad  (m(x),A)=\frac{(\dif\rho_v)_eA}{||v||^2}, \quad 0\ne v\in V,\;x=[v],
\end{equation}
where  $(\cdot,\cdot)$ is an $\Ad(K)$-invariant real inner product on $\im\kg$ and $\PP V$ is the projective
space of lines in $V$.  If $\pi:V\setminus\{ 0\}\mapsto \PP V$ denotes the usual projection map, then
$\pi(v)=[v]=x$. Under the natural identifications $\im\kg=\im\kg^*=\kg^*$, the function $m$ is the moment map
from symplectic geometry, corresponding to the Hamiltonian action of $K$ on the symplectic manifold $\PP V$ (see
for instance the survey \cite{Krw} or \cite[Chapter 8]{MFK}).

\begin{definition}\label{stable} {\rm
A vector $v\in V$ is said to be {\it unstable} if $0\in\overline{G.v}$, the closure of the orbit $G.v$, and {\it
semistable} if $0\notin\overline{G.v}$.  The set $N\subset V$ of unstable vectors is called the {\it null-cone}
and the set of semistable vectors is denoted by $V^{\rm ss}$. }
\end{definition}

The set $V^{\rm ss}$ is open in $V$.  A well known result due to G. Kempf and L. Ness \cite{KmNs} asserts that
an orbit $G.v$ is closed if and only if $G.[v]$ meets $m^{-1}(0)$, and in that case the intersection coincides
with $K.[v]$.  For any $v\in V^{\rm ss}$ there is a unique closed orbit in $\overline{G.v}$, and we then say
that $v\sim w$ if the closures of their orbits contain the same closed orbit.  If $\PP V^{\rm ss}=\pi(V^{\rm
ss})$, then the so called {\it categorical quotient} $\PP V//G=\PP V^{\rm ss}/\sim$ is homeomorphic to
$m^{-1}(0)/K$, the symplectic quotient or reduced (phase) space of the symplectic manifold $\PP V$ at the level
$0$.

Consider the functional square norm of the moment map
\begin{equation}\label{norm}
F:\PP V\mapsto\RR, \qquad  F(x)=||m(x)||^2=(m(x),m(x)).
\end{equation}

Thus an orbit $G.v$ is closed if and only if $F(x)=0$ for some $x\in G.[v]$, and in that case, the set of zeros
of $F|_{G.x}$ coincides with $K.x$.  Moreover, $\PP V^{\rm ss}$ is the set of points $[v]\in\PP V$ with the
property that the limit of the negative gradient flow of $F$ is in $m^{-1}(0)$.  A natural question arises: what
is the role played by the remaining critical points of $F:\PP V\mapsto\RR$ (i.e. those for which $F(x)>0$) in
the study of the $G$-orbit space of the action of $G$ on $\PP V$, as well as on other complex projective
$G$-varieties contained in $\PP V$?.  This was precisely the aim of the paper \cite{Ns} (see also \cite{Krw1}),
where it is shown that the non-minimal critical points have influence in the study of the orbit space of the
null-cone.

\begin{theorem}\cite{Ns}\label{ness1}
Let $F|_{G.x}$ denote the restriction of $F:\PP V\mapsto\RR$ to the $G$-orbit of $x\in \PP V$.  If $x$ is a
critical point of $F$ then
\begin{itemize}
\item[(i)] $F|_{G.x}$ attains its minimum value at $x$.

\item[(ii)] $F|_{G.x}$ attains its minimum value only on the orbit $K.x$.
\end{itemize}
The non-minimal critical points (i.e. $F(x)>0$) are all in the null-cone, that is, $0\in\overline{G.v}$,
$x=[v]$.
\end{theorem}

The proof of part (i) of the above theorem is based in the presently well known convexity properties of moment
maps (see \cite{Krw}).

\begin{theorem}\cite{Ns}\label{ness2}
The negative gradient flow of $F:\PP V\mapsto\RR$ determines a stratification of the null-cone $N$.  A stratum
$S_{\la h\ra}$ of $N$, $ h\in\im\kg$, is the set of all the points $x\in \PP V$ which flow into $C_{\la h\ra}$,
where $C_{\la h\ra}$ is the set of critical points $y$ of $F$ such that $m^*(y)\in\Ad(K). h$.  Moreover,
\begin{itemize}
\item[(i)] For each stratum $S_{\la h\ra}$ there exists a subspace $V_ {\la h\ra}\subset V$ and a reductive
subgroup $G_{\la h\ra}\subset G$ such that   $C_{\la h\ra}/K =\PP V_{\la h\ra}\mum G_{\la h\ra}$, and hence
$C_{\la h\ra}/K$ is a projective algebraic variety.

\item[(ii)] There are finitely many strata.

\item[(iii)] Each stratum $S_{\la h\ra}$ is Zariski-locally closed, irreducible and nonsingular.
\end{itemize}
\end{theorem}

Therefore, one may conclude that there exist also distinguished orbits in the null-cone $N$, namely those
containing a critical point of $F$, which would play in some sense the same role as the closed orbits in $V^{\rm
ss}$.  If $C$ denotes the set of all critical points of the functional $F=||m||^2:\PP V\mapsto\RR$, then it is
natural to define the {\it Kirwan-Ness quotient} for the action of $G$ on a projective $G$-variety $X\subset \PP
V$  by
$$
X\kir G := C_X/K, \qquad C_X=C\cap X.
$$
It follows from Theorem \ref{ness1}, (ii), that $X\kir G=GC_X/G$ and thus $X\kir G$ parameterizes precisely
those $G$-orbits containing a critical point of $F$.  We note that if the action of $G$ on $V$ contains the
homotheties $\{ v\mapsto tv:t\in\CC^*\}$, then we may also define the Kirwan-Ness quotient for the action of $G$
on an algebraic $G$-variety $X\subset V$ by
$$
X\kir G=\{ v\in X:[v]\in C\}/\CC^*K,
$$
and clearly $X\kir G=\pi(X)\kir G$.  This new quotient $X\kir G$ is not a projective algebraic variety (not even
Hausdorff) as in the case of $X\mum G$, but it still has some nice properties.  Let $\bca$ denote the finite set
indexing the strata described in Theorem \ref{ness2}.  Consider the following partial order in $\bca$:
$\alpha<\beta$ if $F([v])<F([w])$ for $[v]\in C_{\alpha}$, $[w]\in C_{\beta}$.  We have that
$$
\overline{S}_{\alpha}\subset S_{\alpha}\cup\bigcup_{\beta>\alpha}S_{\beta},
$$
and therefore the Kirwan-Ness quotient can be decomposed as the disjoint union
$$
X\kir G=\bigcup_{\alpha\in\bca}X_{\alpha}
$$
of the projective algebraic varieties $X_{\alpha}=(C_{\alpha}\cap X)/K$, satisfying the following frontier
property:
$$
\overline{X}_{\alpha}\subset X_{\alpha}\cup\bigcup_{\beta>\alpha}X_{\beta}.
$$
Recall that for $\alpha=\la 0\ra$, we have $X_{\alpha}=\PP V\mum G$, the categorical quotient.

Fix a maximal torus $T\subset G$ and a Borel subgroup $T\subset B\subset G$.  Thus $B$ determines a positive
Weyl chamber $\im\tg^+\subset\im\tg$, where $\tg$ is the Lie algebra of $T$.  Since $m^*:\PP V\mapsto\im\kg$ is
$K$-equivariant, for any $x\in \PP V$, $m^*(K.x)$ has a unique point of intersection $\check{m}(x)$ with
$\im\tg^+$, thus determining a function
$$
\check{m}:\PP V\mapsto\im\tg^+.
$$

\begin{theorem}\cite{Mm}\label{convex}
The image via $\check{m}$ of any $G$-invariant closed set in $\PP V$ is a rational convex polytope in
$\im\tg^+$, that is, the convex hull of a finite number of points with rational coordinates.
\end{theorem}

We also have that $F(x)$ is precisely the square of the distance from the origin in $\im\kg$ to the point
$\check{m}(x)$ in the polytope $\check{m}(\overline{G.x})$.

\section{The moment map for skew-symmetric algebras}\label{bil}

The object of this section is to study the notions described in Section \ref{ness} for the $\Gl(n)$-action on
the vector space where the Lie algebras live.  Some of the results given here follow readily from the general
case proved in \cite{Ns}, but since the adaptation can be sometimes rather difficult, we have argued directly in
our particular situation for completeness.

Let $V_n=\Lambda^2(\CC^n)^*\otimes\CC^n$ be the vector space of all alternating bilinear maps from
$\CC^n\times\CC^n$ to $\CC^n$, or in other words, the space of all skew-symmetric (non-associative) algebras of
dimension $n$.  There is a natural action of $\Gl(n)=\Gl(n,\CC)$ on $V_n$ given by
\begin{equation}\label{action}
g.\mu(X,Y)=g\mu(g^{-1}X,g^{-1}Y), \qquad X,Y\in\CC^n, \quad g\in\Gl(n),\quad \mu\in V_n.
\end{equation}
We note that any $\mu\in V_n$ is in the null-cone, since for $g_t=tI$ we have that
$$
\lim_{t\to 0}g_t^{-1}.\mu=\lim_{t\to 0}t\mu=0,
$$
and so $0\in\overline{\Gl(n).\mu}$.  The usual hermitian inner product $\ip$ on $\CC^n$, defines a
$\U(n)$-invariant hermitian inner product on $V_n$, denoted also by $\ip$, as follows:
\begin{equation}\label{innV}
\la\mu,\lambda\ra=\sum_{ijk}\la\mu(X_i,X_j),X_k\ra\overline{\la\lambda(X_i,X_j),X_k\ra},
\end{equation}
where $\{ X_1,...,X_n\}$ is any orthonormal basis of $\CC^n$.  The Lie algebra of $\Gl(n)$ decomposes
$\glg(n)=\ug(n)+\im\ug(n)$ in skew-hermitian and hermitian transformations respectively, and an
$\Ad(\U(n))$-invariant hermitian inner product on $\glg(n)$ is given by
\begin{equation}\label{inng}
(A,B)=\tr{AB^*}, \qquad A,B\in\glg(n).
\end{equation}
Thus we use $(\cdot,\cdot)$ to identify $\im\ug(n)$ with $\im\ug(n)^*$.  For each $\mu\in V_n$, consider the
hermitian map $\Rin_{\mu}\in\im\ug(n)$ defined by
\begin{equation}\label{defR2}
\Rin_{\mu}=-4\displaystyle{\sum_{i}}(\ad_{\mu}{X_i})^*\ad_{\mu}{X_i}
+2\displaystyle{\sum_{i}}\ad_{\mu}{X_i}(\ad_{\mu}{X_i})^*,
\end{equation}
where the adjoint map $\ad_{\mu}{X}:\CC^n\mapsto\CC^n$ or left multiplication by $X$ of the algebra $\mu$ is
given, as usual, by $\ad_{\mu}{X}(Y)=\mu(X,Y)$.  It is a simple calculation to see that
\begin{equation}\label{defR}
\begin{array}{rl}
\la\Rin_{\mu}X,Y\ra=&-4\displaystyle{\sum_{ij}}\la\mu(X,X_i),X_j\ra\overline{\la\mu(Y,X_i),X_j\ra} \\
&+2\displaystyle{\sum_{ij}}\overline{\la\mu(X_i,X_j),X\ra}\la\mu(X_i,X_j), Y\ra,
\end{array}
\end{equation}
for all $X,Y\in\CC^n$.  We will see below that the map $\mu\mapsto\Rin_{\mu}$ is precisely the moment map for
the action (\ref{action}).  The action of $\glg(n)$ on $V_n$ obtained by differentiation of (\ref{action}) is
given by
\begin{equation}\label{actiong}
A.\mu=-\delta_{\mu}(A)= A\mu(\cdot,\cdot)-\mu(A\cdot,\cdot)-\mu(\cdot,A\cdot), \qquad A\in\glg(n),\quad\mu\in
V_n.
\end{equation}
If $\mu\in V_n$ satisfies the Jacobi condition, then $\delta_{\mu}:\glg(n)\mapsto V_n$ coincides with the
cohomology coboundary operator of the Lie algebra $(\CC^n,\mu)$ from level $1$ to $2$, relative to cohomology
with values in the adjoint representation.  Recall that $\Ker{\delta_{\mu}}=\Der(\mu)$, the Lie algebra of
derivations of the algebra $\mu$.  The following technical lemma will be crucial in several proofs throughout
the paper.

\begin{lemma}\label{dR}
Let $p:\glg(n)\mapsto\im\ug(n)$ be the projection relative to the decomposition $\glg(n)=\ug(n)+\im\ug(n)$, and
let $\delta_{\mu}^*:V_n\mapsto\glg(n)$ denote the transpose of $\delta_{\mu}:\glg(n)\mapsto V_n$, relative to
the hermitian inner products $(\cdot,\cdot)$ and $\ip$, on $\glg(n)$ and $V_n$ respectively.
\begin{itemize}
\item[(i)] If $\Rin:V_n\mapsto\im\ug(n)$ is defined by $\Rin(\mu)=\Rin_{\mu}$ for every $\mu\in V_n$, then
$$
(\dif\Rin)_{\mu}=-4p\circ\delta_{\mu}^*.
$$
\item[(ii)] $\tr{\Rin_{\mu}D}=0$ for any $D\in\im\ug(n)\cap\Der(\mu)$.

\item[(iii)] $\tr{\Rin_{\mu}[A,A^*]}\geq 0$ for any $A\in\Der(\mu)$.  Equality holds if and only if we also have
that $A^*\in\Der(\mu)$.
\end{itemize}
\end{lemma}

\begin{proof} (i).  Consider the line $\mu+t\lambda$ with $\mu,\lambda\in V_n$, $t\in\RR$.  Using (\ref{defR}), for
any $A\in\im\ug(n)$ we obtain that
$$
\begin{array}{rl}
((\dif{\Rin})_{\mu}(\lambda),A) =&
\tr{A(\dif{\Rin})_{\mu}(\lambda)}= \displaystyle{\sum_{pr}}\la (\dif{\Rin})_{\mu}(\lambda)X_p,X_r\ra\overline{\la AX_p,X_r\ra} \\

=&\displaystyle{\sum_{pr}}\ddt_{|_0}\la\Rin_{\mu+t\lambda}X_p,X_r\ra\overline{\la AX_p,X_r\ra} \\

=&\displaystyle{\sum_{pr}}\Big(\displaystyle{\sum_{ij}} -4\la\lambda(X_p,X_i),X_j\ra
\overline{\la\mu(X_r,X_i),X_j\ra} \\

&-4\la\mu(X_p,X_i),X_j\ra
\overline{\la\lambda(X_r,X_i),X_j\ra} \\ \\

&+2\overline{\la\lambda(X_i,X_j),X_p\ra}\la\mu(X_i,X_j),X_r\ra\\ \\

&+2\overline{\la\mu(X_i,X_j),X_p\ra}\la\lambda(X_i,X_j),X_r\ra\Big)\overline{\la AX_p,X_r\ra}.
\end{array}
$$
We now interchange, in even lines, the indeces $p$ and $r$, obtaining that
$$
\begin{array}{rl}
((\dif{\Rin})_{\mu}(\lambda),A) =& \Re\displaystyle{\sum_{prij}}-8\la\overline{\mu(X_r,X_i),X_j\ra}
\la\lambda(X_p,X_i),X_j\ra\overline{\la AX_p,X_r\ra}\\

&+4\overline{\la\mu(X_i,X_j),X_r\ra}\la\lambda(X_i,X_j),X_p\ra\la AX_p,X_r\ra.
\end{array}
$$
On the other hand, we have that
$$
\begin{array}{rl}
\la\lambda,\delta_{\mu}(A)\ra

=&\displaystyle{\sum_{pij}}
\la\lambda(X_p,X_i),X_j\ra\overline{\la\delta_{\mu}(A)(X_p,X_i),X_j\ra} \\

=&\displaystyle{\sum_{pij}}\la\lambda(X_p,X_i),X_j\ra\overline{\la\mu(AX_p,X_i),X_j\ra} \\

&+ \la\lambda(X_p,X_i),X_j\ra\overline{\la\mu(X_p,AX_i),X_j\ra}  \\ \\

&-\la\lambda(X_p,X_i),X_j\ra\overline{\la A\mu(X_p,X_i),X_j\ra} \\ \\

=&
\displaystyle{\sum_{pijr}}\la\lambda(X_p,X_i),X_j\ra\overline{\la\mu(X_r,X_i),X_j\ra}\overline{\la AX_p,X_r\ra} \\

& +\la\lambda(X_p,X_i),X_j\ra\overline{\la\mu(X_p,X_r),X_j\ra}\overline{\la AX_i,X_r\ra}  \\ \\

& -\la\lambda(X_p,X_i),X_j\ra\overline{\la\mu(X_p,X_i),X_r\ra}\la AX_j,X_r\ra.
\end{array}
$$
By interchanging the indexes $p$ and $i$ in the second line, and $p$ and $j$ in the third one, we get
$$
\begin{array}{rl}
\la\lambda,\delta_{\mu}(A)\ra =& 2\displaystyle{\sum_{prij}}\la\lambda(X_p,X_i),X_j\ra
\la\overline{\mu(X_r,X_i),X_j\ra}
\overline{\la AX_p,X_r\ra}\\

&-\la\lambda(X_i,X_j),X_p\ra\overline{\la\mu(X_i,X_j),X_r\ra}\la AX_p,X_r\ra.
\end{array}
$$
We then can deduce from the two computations above that
$$
((\dif{\Rin})_{\mu}(\lambda),A)=-4\Re\la\lambda,\delta_{\mu}(A)\ra= -4\Re(\delta_{\mu}^*(\lambda),A)=
\Re(-4p\circ\delta_{\mu}^*(\lambda),A),
$$
for every $A\in\im\ug(n)$, which concludes the proof of (i).

\no (ii).  Using that $2\Rin_{\mu}=\ddt_{|_0}(1+t)^2\Rin_{\mu}=\ddt|_0\Rin_{\mu+t\mu}=(\dif\Rin)_{\mu}(\mu)$ and
part (i), we obtain that
$$
\begin{array}{rl}
\tr{\Rin_{\mu}D}=&\unm\tr{(\dif\Rin)_{\mu}(\mu)D}=-2\tr{p\circ\delta_{\mu}^*(\mu)D}\\ \\
=&-2\Re(\delta_{\mu}^*(\mu),D)=-2\Re\la\mu,\delta_{\mu}(D)\ra=0.
\end{array}
$$

\noindent (iii).  It follows easily from the $K$-invariance of $\ip$ on $V_n$ that
$$
\la A.\mu,\lambda\ra=\la\mu,A^*.\lambda\ra
$$
for any $A\in\glg(n)$.  We then have by part (i) that
$$
\begin{array}{rl}
\tr{\Rin_{\mu}[A,A^*]} &= \unm\tr{(\dif\Rin)_{\mu}(\mu)[A,A^*]}=
-2\Re(\delta_{\mu}^*(\mu),[A,A^*]) \\ \\
&=-2\Re\la\mu,\delta_{\mu}([A,A^*])\ra
= 2\Re\la\mu,[A,A^*].\mu\ra \\ \\
&=2\Re\la\mu,A.(A^*.\mu)-A^*.(A.\mu)\ra
=2\Re\la\mu,A.(A^*.\mu)\ra \\ \\

&=2\la A^*.\mu,A^*.\mu\ra\geq 0,
\end{array}
$$
and it equals $0$ if and only if $A^*.\mu=0$, that is, $A^*\in\Der(\mu)$.
\end{proof}

We now calculate the moment map and the gradient of the functional square norm of the moment map, obtaining a
rather computable characterization of their critical points, which will be very useful.

\begin{proposition}\label{calculo}
The moment map $m:\PP V_n\mapsto\im\ug(n)$, the functional square norm of the moment map $F_n=||m||^2:\PP
V_n\mapsto\RR$ and the gradient of $F$ are respectively given by
\begin{equation}\label{grad}
m([\mu])=\Rin_{\mu},  \quad F_n([\mu])=\tr{\Rin_{\mu}^2},  \quad
\grad(F_n)_{[\mu]}=-8\pi_{\ast}\delta_{\mu}(\Rin_{\mu}), \qquad ||\mu||=1,
\end{equation}
where $\pi_{\ast}$ denotes the derivative of $\pi:V_n\setminus\{ 0\}\mapsto \PP V_n$, the canonical projection.
Moreover, the following statements are equivalent:
\begin{itemize}
\item[(i)] $[\mu]\in \PP V_n$ is a critical point of $F_n$. \item[(ii)] $[\mu]\in \PP V_n$ is a critical point
of $F_n|_{\Gl(n).[\mu]}$. \item[(iii)] $\Rin_{\mu}=cI+D$ for some $c\in\RR$ and $D\in\Der(\mu)$.
\end{itemize}
\end{proposition}

\begin{proof} For any $A\in\im\ug(n)$ we have that
$$
(\dif\rho_{\mu})_I(A)=\ddt|_0\la e^{tA}.\mu,e^{tA}.\mu\ra =-2\Re\la\delta_{\mu}(A),\mu\ra
=-2\Re(A,\delta_{\mu}^*(\mu)).
$$
Now, using Lemma \ref{dR},(i), we get
$$
(\dif\rho_{\mu})_I(A)=-2\Re(A,-\unc(\dif\Rin)_{\mu}(\mu)\ra =\Re(A,\Rin_{\mu}) =\Re\tr{\Rin_{\mu}A},
$$
which proves the formula for $m$.  The first assertion on $F_n$ is self-evident.  To prove the second one, we
only need to compute the gradient of $F_n:V_n\mapsto\RR$, $F_n(\mu)=\tr{\Rin_{\mu}^2}$, and then to project it
via $\pi_{\ast}$.  If $\mu,\lambda\in V_n$, then it follows from Lemma \ref{dR}, (i), that
$$
\begin{array}{rl}
\Re\la\grad(F_n)_{\mu},\lambda\ra & =\ddt|_{0}F_n(\mu+t\lambda)

=\ddt|_0\tr{\Rin_{\mu+t\lambda}^2}

=2\Re\tr{(\ddt|_0\Rin_{\mu+t\lambda})\Rin_{\mu}} \\ \\

&=2\Re((\dif\Rin_{\mu})_{\mu}(\lambda),\Rin_{\mu})

=-8\Re(\delta_{\mu}^*(\lambda),\Rin_{\mu})

=-8\Re\la\lambda,\delta_{\mu}(\Rin_{\mu})\ra.
\end{array}
$$
This implies that $\grad(F_n)_{\mu}=-8\delta_{\mu}(\Rin_{\mu})$, concluding the proof of (\ref{grad}).

The equivalence between (i) and (ii) has been observed in \cite{Ns}, and it follows from the fact that
$\grad(F_n)_{[\mu]}\in\tang_{[\mu]}\Gl(n).[\mu]$ for any $[\mu]\in \PP V_n$.  Indeed, by (\ref{grad}) we have
that
$$
\grad(F_n)_{[\mu]}=8\pi_{\ast}(\ddt|_0e^{t\Rin_{\mu}}.\mu)= 8\ddt|_0e^{t\Rin_{\mu}}.[\mu].
$$
In order to prove that (i) is equivalent to (iii), recall first that $\Ker{\pi_{\ast}}(\mu)=\CC\mu$.  Hence we
obtain from (\ref{grad}) that $[\mu]$ is a critical point of $F_n$ if and only if
$\delta_{\mu}(\Rin_{\mu})\in\CC\mu$, or equivalently, $\Rin_{\mu}\in\CC I\oplus\Der(\mu)$, since
$\delta_{\mu}(I)=\mu$ and $\Ker{\delta_{\mu}}=\Der(\mu)$ (see (\ref{actiong})).  If $\Rin_{\mu}=cI+D$ with
$c\in\CC$ and $D\in\Der(\mu)$, then it is evident that $D$ is normal, and so $D^*$ is also a derivation of $\mu$
(see Lemma \ref{dR}, (iii)).  Let $D=D_s+D_h$ be the skew-hermitian and hermitian parts of $D$.  Since
$D_s=-\im\Im(c)I$ has to be a derivation we get that $D_s=0$, and therefore $\Rin_{\mu}=cI+D_h$ and $c\in\RR$,
as it was to be shown.
\end{proof}

In the frame of skew-symmetric algebras, the result due to L. Ness given in Theorem \ref{ness1} can be stated as
follows.

\begin{theorem}\label{min}\cite{Ns}
If $[\mu]$ is a critical point of the functional $F_n:\PP V_n\mapsto\RR$ given by
$F_n([\lambda])=\tr{\Rin_{\lambda}^2}$ ($||\lambda||=1$), then
\begin{itemize}
\item[(i)] $F_n|_{\Gl(n).[\mu]}$ attains its minimum value at $[\mu]$.

\item[(ii)] $[\lambda]\in \Gl(n).[\mu]$ is a critical point of $F_n$ if and only if $[\lambda]\in\U(n).[\mu]$.
\end{itemize}
\end{theorem}

We now describe some particular features of the critical points and the stratification for our case.

\begin{lemma}\label{cmu}
Let $[\mu]\in \PP V_n$ be a critical point of $F_n$, say $\Rin_{\mu}=c_{\mu}I+D_{\mu}$ for some $c_{\mu}\in\RR$
and $D_{\mu}\in\Der(\mu)$.  Then
$c_{\mu}=\frac{\tr{\Rin_{\mu}^2}}{\tr{\Rin_{\mu}}}=-\unm\frac{\tr{\Rin_{\mu}^2}}{||\mu||^2}$, and if $D_{\mu}\ne
0$ then $\tr{D_{\mu}}>0$ and $c_{\mu}=-\frac{\tr{D_{\mu}^2}}{\tr{D_{\mu}}}$.
\end{lemma}

\begin{proof}
Both assertions follow from the fact that $\tr{\Rin_{\mu}D_{\mu}}=0$ (see Lemma \ref{dR},(ii)) and
$\tr{\Rin_{\mu}}=-2||\mu||^2$.
\end{proof}

The proof of the following rationality result, which is just a bit stronger than the given in \cite[Section
4]{Ns} for the general case, is based on the proof by J. Heber of \cite[Thm 4.14]{H}.

\begin{theorem}\label{rationalization}
Let $[\mu]\in \PP V_n$ be a critical point of $F_n$, with $\Rin_{\mu}=c_{\mu}I+D_{\mu}$ for some $c_{\mu}\in\RR$
and $D_{\mu}\in\Der(\mu)$. Then there exists $c>0$ such that the eigenvalues of $cD_{\mu}$ are all nonnegative
integers prime to each other, say $k_1<...<k_r\in\ZZ_{\geq 0}$ with multiplicities $d_1,...,d_r\in\NN$.
\end{theorem}

\begin{proof} If $D_{\mu}=0$ then there is nothing to prove, so we assume $D_{\mu}\ne 0$.  Consider $\CC^n=\ggo_1\oplus...\oplus\ggo_r$ the orthogonal decomposition in eigenspaces of $D_{\mu}$, say $D_{\mu}|_{\ggo_i}=c_iI_{\ggo_i}$ with
$c_1<...<c_r$.  Since $D_{\mu}\in\Der(\mu)$ we have that $\mu(\ggo_i,\ggo_j)\subset\ggo_k$ if and only if
$c_i+c_j=c_k$; otherwise $\mu(\ggo_i,\ggo_j)=0$.  A crucial point here is that any hermitian map $A$ of $\CC^n$
defined by $A|_{\ggo_i}=a_iI_{\ggo_i}$, satisfying $a_i+a_j=a_k$ for all $i,j,k$ such that $c_i+c_j=c_k$, is
also a derivation of $\mu$.  Now, using Lemma \ref{dR},(ii) and Lemma \ref{cmu}, we get
\begin{equation}\label{eigen1}
\tr{D_{\mu}A}=c_{\mu}\tr{A}=-\frac{\tr{D_{\mu}^2}}{\tr{D_{\mu}}}\tr{A}
\end{equation}
for every $A$ satisfying the above conditions.  In other words, if $\{ e_1,...,e_r\}$ denotes the canonical
basis of $\RR^r$ and $\alpha=\sum c_i^2/\sum c_i$, then (\ref{eigen1}) says that the vector
$(c_1-\alpha,...,c_r-\alpha)$ is orthogonal to $F^{\perp}$, where
$$
F=\{ e_i+e_j-e_k:c_i+c_j=c_k\}.
$$
This implies that $(c_1-\alpha,...,c_r-\alpha)\in\la F\ra$ (subspace of $\RR^r$ linearly spanned by $F$), and
thus
$$
(c_1-\alpha,...,c_r-\alpha)=\sum_{p=1}^sb_p(e_{i_p}+e_{j_p}-e_{k_p}), \qquad b_p\in\RR,
$$
where $\{ e_{i_p}+e_{j_p}-e_{k_p} : p=1,...,s\}$ is a basis of $\la F\ra$.  We now consider the $(s\times
r)$-matrix
$$
E=\left[\begin{array}{c}
e_{i_1}+e_{j_1}-e_{k_1} \\
\vdots \\
e_{i_s}+e_{j_s}-e_{k_s}
\end{array}\right], \; \mbox{so that}\; EE^t\in\Gl(s,\QQ),
$$
and furthermore
$$
\quad E \left(\begin{array}{c}
c_1 \\
\vdots \\
c_r
\end{array}\right)=0, \quad E \left(\begin{array}{c}
1 \\
\vdots \\
1
\end{array}\right)= \left(\begin{array}{c}
1 \\
\vdots \\
1
\end{array}\right) \; \mbox{and}\; E^t\left(\begin{array}{c}
b_1 \\
\vdots \\
b_s
\end{array}\right)=(c_1-\alpha,...,c_r-\alpha).
$$
This implies that
\begin{equation}\label{eigen2}
\frac{1}{\alpha}(c_1,...,c_r)=(1,...,1)-E^t(EE^t)^{-1}\left(\begin{array}{c}
1 \\
\vdots \\
1
\end{array}\right)\in\QQ^r,
\end{equation}
as it was to be shown.  It remains to prove that $0\leq c_1$.  If $D_{\mu}X=c_1X$ then
$c_1\ad_{\mu}{X}=\ad_{\mu}{D_{\mu}X}=[D_{\mu},\ad_{\mu}{X}]$, and therefore
$$
c_1\tr{\ad_{\mu}{X}(\ad_{\mu}{X})^*}=\tr{D_{\mu}[\ad_{\mu}{X},(\ad_{\mu}{X})^*]}=\tr{\Rin_{\mu}[\ad_{\mu}{X},(\ad_{\mu}{X})^*]}\geq
0,
$$
by Lemma \ref{dR},(iii).  If we assume $\ad_{\mu}{X}\ne 0$, we obtain that $c_1\geq 0$ .  Otherwise, if
$\ad_{\mu}{X}=0$ then $0\leq\la\Rin_{\mu}X,X\ra=(c_{\mu}+c_1)\la X,X\ra$, and so $c_1>0$ since $c_{\mu}<0$ by
Lemma \ref{cmu}.
\end{proof}

\begin{definition}\label{tipo}
{\rm The data set $\eigen$ in the above theorem is called the {\it type} of the critical point $[\mu]$.  }
\end{definition}

\begin{proposition}\label{finitos}
In any fixed dimension $n$, there are only finitely many types of critical points of $F_n:\PP V_n\mapsto\RR$ .
\end{proposition}

\begin{proof}
Equation (\ref{eigen2}) says that the numbers $k_1,...,k_r$ can be recovered from the knowledge of $r$ and the
finite set $F_n$ alone.  Thus the finiteness of the types follows from the finiteness of the partitions
$n=d_1+...+d_r$ and the different sets $F_n$.
\end{proof}

Since the null-cone of $V_n$ is all of $V_n$, the strata $S_{\la h\ra}$'s described in Theorem \ref{ness2}
determines a stratification of $V_n$.  The set of types of critical points is in bijection with the set of
strata.  Indeed, for a type $ \alpha=\eigen$ we define
$$
h_{\alpha}=-\frac{k_1^2d_1+...+k_r^2d_r}{k_1d_1+...+k_rd_r}I+ \left[\begin{array}{ccc}
k_1I_{d_1}&&\\
&\ddots&\\
&&k_rI_{d_r}
\end{array}\right]\in\im\ug(n),
$$
where $I_{d_i}$ denotes the $d_i\times d_i$ identity matrix.  Thus the set $C_{ \alpha}=C_{\la h_{\alpha}\ra}$
given in Theorem \ref{ness2} is precisely the set of critical points $[\mu]$ of $F_n:\PP V_n\mapsto\RR$ such
that $m([\mu])=\Rin_{\mu}=c_{\mu}I+D_{\mu}$ is conjugate to $h_{\alpha}$, or equivalently,
$$
C_{ \alpha}=\{ [\mu]\in\PP V_n:[\mu]\;\mbox{is a critical point of}\; F_n\;\mbox{of type}\;  \alpha\}.
$$
This implies that each stratum is of the form $S_{ \alpha}=S_{\la h_{\alpha}\ra}$ for some type $ \alpha$, where
\begin{equation}\label{stratum}
S_{ \alpha}=\left\{\begin{array}{l} [\lambda]\in\PP V_n:\;\mbox{the limit of the}\;-\grad(F_n)\\ \mbox{flow
starting from}\; [\lambda]\;\mbox{is in}\; C_{ \alpha}
\end{array}\right\}.
\end{equation}
If $\alpha_1,...,\alpha_r$ are the different types of critical points then
$$
\PP V_n=S_{\alpha_1}\cup...\cup S_{\alpha_s}
$$
determines a stratification of $\PP V_n$, for which each stratum $S_{\alpha_i}$ is locally closed, irreducible
and nonsingular.  If $\check{m}:\PP V_n\mapsto \im\ug(n)^+$ is the function considered in Theorem \ref{convex},
then $\check{m}([\mu])=h_{\alpha}$ for any critical point $[\mu]$ of type $\alpha$.

Since $\grad(F_n)$ is always tangent to $\Gl(n)$-orbits, we have that the stratification of any
$\Gl(n)$-invariant projective algebraic variety $X\subset\PP V_n$ is obtained just by intersecting the stratum
$S_{ \alpha}$ of $\PP V_n$ with $X$.  In view of the equivalence between (i) and (ii) in Proposition
\ref{calculo}, the critical points of $F_n:X\mapsto\RR$ are precisely the critical points of $F_n:\PP
V_n\mapsto\RR$ which lie in $X$.  Thus we will often refer to them just as critical points of $F_n$.

The following two propositions on the critical values of $F_n$ and abelian factors were proved in \cite[Section
3]{L1}.

\begin{proposition}\label{values}
Let $[\mu]\in\PP V_n$ be a critical point of $F_n$ of type $\eigen$, different from $(0,n)$.  Then,
$$
F_n([\mu])=4\Big(n-\frac{(k_1d_1+...+k_rd_r)^2}{k_1^2d_1+...+k_r^2d_r}\Big)^{-1}.
$$
\end{proposition}

\begin{proposition}\label{factor}
Let $[\mu]\in\PP V_n$ and $[\lambda]\in\PP V_m$ be critical points of $F_n$ and $F_m$ respectively.  Then the
direct sum $[\mu\oplus c\lambda]\in\PP V_{n+m}$ is a critical point of $F_{n+m}$ for a suitable $c\in\RR$.  If
in addition $\lambda$ is abelian, then $F_n([\mu])=F_{n+m}([\mu\oplus\lambda])$ and the type of
$[\mu\oplus\lambda]$ is given by
$$
(ak_1<...<\frac{k_1^2d_1+...+k_r^2d_r}{d}<...<ak_r;d_1,...,m,...,d_r),
$$
where $d={\rm mcd}(k_1d_1+...+k_rd_r,k_1^2d_1+...+k_r^2d_r)$ and $a=\frac{k_1d_1+...+k_rd_r}{d}$.  In case that
$\frac{k_1^2d_1+...+k_r^2d_r}{d}=ak_i$ for some $i$, then the multiplicity is $m+d_i$.
\end{proposition}

\section{Critical points in the variety of Lie algebras}\label{critvar}

The space of all $n$-dimensional complex Lie algebras can be naturally identified with the subset $\lca_n\subset
V_n$ of all Lie brackets.  $\lca_n$ is an algebraic set, since the Jacobi identity is given by polynomial
conditions.  The isomorphism class of a Lie algebra $\mu\in\lca_n$ is then the orbit $\Gl(n).\mu$ under the
`change of basis' action of $\Gl(n)$ on $\lca_n$ given in (\ref{action}).

\begin{definition}\label{deg} {\rm
We say that $\mu$ {\it degenerates to} $\lambda$ if $\lambda\in\overline{\Gl(n).\mu}$, which will be often
denoted by $\mu\mapsto\lambda$.  }
\end{definition}

Every degeneration will be assumed to be nontrivial, that is, $\lambda$ lies in the boundary of $\Gl(n).\mu$.
If $\mu\mapsto\lambda$ then we may say roughly that $\lambda$ is `more abelian' than $\mu$; in fact,
$\dim{\Der(\lambda)}>\dim{\Der(\mu)}$, $\dim{\lambda(\CC^n,\CC^n)}\leq\dim{\mu(\CC^n,\CC^n)}$,
$\dim{\zg(\lambda)}\geq\dim{\zg(\mu)}$ and $\abel(\lambda)\geq\abel(\mu)$, where $\zg(\mu)$ denotes the center
of $\mu$ and $\abel(\mu)$ is the dimension of a maximal abelian subalgebra of $\mu$ (see \cite{BS}).

We note that the $-\grad(F_n)$ flow $\mu(t)$ defines a degeneration (possibly trivial) of the starting point
$\mu(0)$ to a critical point of $F_n$, since $\mu(t)\in\Gl(n).\mu(0)$ for any $t$.  Therefore, such a
distinguished degeneration associates to each Lie algebra $\mu=\mu(0)$ a $\ZZ_{\geq 0}$-graded Lie algebra
$\lambda=\lim_{t\to 0}\mu(t)$, and there are finitely many possible gradations.  It is easy to see that if $\mu$
is nilpotent then the associated $\lambda$ is actually $\NN$-graded.

Let $L_n$ be the projective algebraic variety obtained by projectivization of the algebraic variety
$\lca_n\subset V_n$ of $n$-dimensional Lie algebras, that is, $L_n=\pi(\lca_n\setminus\{ 0\})$.  We are
intrigued by the Kirwan-Ness quotient $L_n\kir\Gl(n)$, and thus we want to understand the critical points of
$F_n$ which lie in $L_n$, where $F_n$ is the functional square norm of the moment map (see $\S$\ref{bil}).  A
first natural question would be to find the global minima and maxima of $F_n:L_n\mapsto\RR$.

\begin{lemma}\label{mini}
$[\lambda]$ is a critical point of $F_n$ of type $(0;n)$ (i.e. $\Rin_{\lambda}\in\RR I$) if and only if
$F_n([\lambda])=\frac{4}{n}$. In that case, $F_n:\PP V_n\mapsto\RR$ attains its minimum at $[\lambda]$.
\end{lemma}

\begin{proof}
It follows from the formula
$$
F_n([\mu])=4\frac{\tr{\Rin_{\mu}^2}}{(\tr{\Rin_{\mu}})^2}, \qquad [\mu]\in\PP V_n,
$$
and a standard analysis of the function $f:\RR^n\mapsto\RR$, $f(c_1,...,c_n)=c_1^2+...+c_n^2/(c_1+...+c_n)^2$.
\end{proof}

\begin{theorem}\label{ss} Assume there exists a semisimple Lie algebra of dimension $n$.  Then $F_n:L_n\mapsto\RR$ attains its minimum value at a point $[\lambda]\in\Gl(n).[\mu]$ if and only if $\mu$ is semisimple.  In such a case, $F_n([\lambda])=\frac{4}{n}$.
\end{theorem}

\begin{proof}
If $\mu$ is semisimple, then by standard methods we can get a basis $\{ Y_i\}$ of $\CC^n$ which is orthonormal
with respect to a suitable hermitian inner product $\ip'$, such that the `Casimir' map satisfies $\sum
(\ad_{\mu}{Y_i})^2\in\RR I$ and $\ad_{\mu}{Y_i}$ is skew-hermitian relative to $\ip'$ for all $i$.  Define
$g\in\Gl(n)$ by $gY_i=X_i$ and consider $\lambda=g.\mu$, where $\{ X_i\}$ is any orthonormal basis of $\CC^n$
relative to our fixed inner product $\ip$.  We note that for any transformation $T$ of $\CC^n$, the adjoints of
$T$ relative to $\ip$ and $\ip'$ are related by $(gTg^{-1})^*=gT^{*'}g^{-1}$.  Thus the moment map at $\lambda$
is given by
$$
\begin{array}{rl}
\Rin_{\lambda}=&
-4\displaystyle{\sum_{i}}(\ad_{\lambda}{X_i})^*\ad_{\lambda}{X_i} +2\displaystyle{\sum_{i}}\ad_{\lambda}{X_i}(\ad_{\lambda}{X_i})^* \\

=&-4\displaystyle{\sum_{i}}(g\ad_{\mu}(g^{-1}X_i)g^{-1})^*g\ad_{\mu}(g^{-1}X_i)g^{-1} \\

& +2\displaystyle{\sum_{i}}g\ad_{\mu}(g^{-1}X_i)g^{-1}(g\ad_{\mu}(g^{-1}X_i)g^{-1})^*  \\

=&-4\displaystyle{\sum_{i}}g(\ad_{\mu}{Y_i})^{*'}g^{-1}g\ad_{\mu}{Y_i}g^{-1} +2\displaystyle{\sum_{i}}g\ad_{\mu}{Y_i}g^{-1}g(\ad_{\mu}{Y_i})^{*'}g^{-1}  \\

=&2g\displaystyle{\sum_{i}}(\ad_{\mu}{Y_i})^2g^{-1}\in\RR I.
\end{array}
$$
This means that $[\lambda]=[g.\mu]\in\Gl(n).[\mu]$ is a critical point of $F_n$ of type $(0;n)$, and thus
$[\lambda]$ is a minimum of $F_n$ (see Lemma \ref{mini}).

Conversely, assume $[\lambda]\in\Gl(n).[\mu]$ is a minimum of $F_n:L_n\mapsto\RR$.  Thus $\Rin_{\lambda}=cI$ for
some $c<0$ (see Lemma \ref{mini}).  It is clear that $\mu$ is semisimple if and only if $\lambda$ is so.  Let
$\sg_{\lambda}$ denote the radical of $\lambda$ and let $\ggo=\hg\oplus\sg_{\lambda}$ the orthogonal
decomposition.  We also decompose orthogonally
$\sg_{\lambda}=\ag\oplus\ngo_{\lambda}=\ag\oplus\vg\oplus\zg_{\lambda}$, where
$\ngo_{\lambda}=\lambda(\sg_{\lambda},\sg_{\lambda})$ and $\zg_{\lambda}$ is the center of the nilpotent Lie
algebra $(\ngo_{\lambda},\lambda|_{\ngo_{\lambda}\times\ngo_{\lambda}})$.  Suppose that $\zg_{\lambda}\ne 0$ and
consider an orthonormal basis $\{ H_i\}$, $\{ A_i\}$, $\{ V_i\}$ and $\{ Z_i\}$ of $\hg$, $\ag$, $\vg$,
$\zg_{\lambda}$ respectively.  If $\{ X_i\}= \{ H_i\}\cup\{ A_i\}\cup\{ V_i\}\cup\{ Z_i\}$ then for every
$Z\in\zg_{\lambda}$,
$$
\begin{array}{rl}
0>\la\Rin_{\lambda} Z,Z\ra = & -4\displaystyle{\sum_{ij}}|\la\lambda(Z,X_i),X_j\ra|^2
+2\displaystyle{\sum_{ij}}|\la\lambda(X_i,X_j),Z\ra|^2 \\
= & -4\displaystyle{\sum_{ij}}|\la\lambda(Z,H_i),Z_j\ra|^2
-4\displaystyle{\sum_{ij}}|\la\lambda(Z,A_i),Z_j\ra|^2 \\
& + 4\displaystyle{\sum_{ij}}|\la\lambda(Z_i,H_j),Z\ra|^2 +4\displaystyle{\sum_{ij}}|\la\lambda(Z_i,A_j),Z\ra|^2
+ \beta(Z),
\end{array}
$$
where $\beta(Z)=2\displaystyle{\sum_{ij}}|\la\lambda(X_i,X_j),Z\ra|^2\geq 0$, $X_i,X_j\in \{ H_i\}\cup\{
A_i\}\cup\{ V_i\}$ (note that both $\ad_{\lambda}\hg$ and $\ad_{\lambda}\ag$ leave invariant  $\zg_{\lambda}$).
This implies that
$$
0> \displaystyle{\sum_{k}}\la\Rin_{\lambda}Z_k,Z_k\ra =\displaystyle{\sum_{k}}\beta(Z_k)\geq 0,
$$
which is a contradiction.  Thus $\zg_{\lambda}$ has to be $\{ 0\}$ and hence $\sg_{\lambda}$ is abelian.  Now,
by applying the same argument to a non-zero $A\in\sg_{\lambda}$ we also get a contradiction, which implies that
$\sg_{\lambda}=0$.  Therefore $\lambda$ is semisimple.
\end{proof}

\begin{remark}\label{sson} {\rm The second part of the above proof implies that any $\mu\in\lca_n$ for which there exists $\lambda\in\Gl(n).\mu$ such that all eigenvalues
of $\Rin_{\lambda}$ are negative, must be semisimple.  In particular, $\Gl(n).\mu$ contains a critical point of
type $(0;n)$ if and only if $\mu$ is semisimple.  Moreover, the stratum $S_{(0;n)}\cap  L_n$ is precisely the
set of semisimple Lie brackets.  Indeed, for any $[\lambda]\in S_{(0;n)}$, the $-\grad(F_n)$ flow
$\{\lambda(t)\}\subset\Gl(n).[\lambda]$ starting from $[\lambda]$ converges to a critical point of type $(0;n)$
and so the eigenvalues of $\Rin_{\lambda(t)}$ will be negative for sufficiently large $t$.         }
\end{remark}

\begin{proposition}\label{red} Assume there does not exists a semisimple Lie algebra of dimension $n$.
\begin{itemize}
\item[(i)] If there exists an $(n-1)$-dimensional semisimple Lie algebra, then $F_n:L_n\mapsto\RR$ attains its
minimum value at some point in $\Gl(n).[\mu]$ if and only if $\mu$ is isomorphic to a reductive Lie algebra with
one-dimensional center.

\item[(ii)] If there is no any $(n-1)$-dimensional semisimple Lie algebra, then $F_n:L_n\mapsto\RR$ attains its
minimum value at some point in $\Gl(n).[\mu]$ if and only if $\mu$ is isomorphic to the direct sum of an
$(n-2)$-dimensional semisimple Lie algebra and the 2-dimensional solvable Lie algebra.
\end{itemize}
In both cases, the type is $(0<1;n-1,1)$ and the minimum value equals $\frac{4}{n-1}$.
\end{proposition}

\begin{proof}
It follows from Proposition \ref{factor} that the orbits considered in both cases contain a critical point of
$F_n$ of type $(0<1;n-1,1)$.  Thus the critical value is $\frac{4}{n-1}$ by Proposition \ref{values}.  It then
suffices to prove that $F_n([\mu])>\frac{4}{n-1}$ for any critical point $[\mu]\in L_n$ of type $\eigen$
different from $(0;n)$, as there is no any semisimple Lie algebra of dimension $n$ (see Remark \ref{sson}).
Actually, we will not use the fact that $\mu\in\lca_n$.

We have that $k_1d_1+...+k_rd_r>k_r$, therefore
$$
\frac{k_1d_1}{\sum k_id_i}k_1+...+\frac{k_rd_r}{\sum k_id_i}k_r\leq k_r<k_1d_1+...+k_rd_r.
$$
This implies that $k_1^2d_1+...+k_r^2d_r<(k_1d_1+...+k_rd_r)^2$, and hence we obtain from Proposition
\ref{values} that
$$
F_n([\mu])=4\left(n-\frac{(k_1d_1+...+k_rd_r)^2}{k_1^2d_1+...+k_r^2d_r}\right)^{-1}<\frac{4}{n-1},
$$
concluding the proof.
\end{proof}

If $\mu$ is semisimple then the minimum $[\lambda]\in\Gl(n).[\mu]$ of $F_n:L_n\mapsto\RR$ is also a minimum of
$F_n:\PP V_n\mapsto\RR$ (see Lemma \ref{mini}).  On the contrary, the minima given in Proposition \ref{red} are
not minima of $F_n:\PP V_n\mapsto\RR$.  In fact, if $[\lambda]\in\PP\lca_{n-z}$ ($z=1$ or $2$) is a semisimple
critical point such that $\Rin_{\lambda}=-12(n-z)I$, then it is easy to check that the bilinear form $\mu\in
V_n$, where we consider $\CC^n=\CC^z\oplus\CC^{n-z}$, defined by
$$
\begin{array}{l}
\mu(X_1,X)=X \quad \forall X\in\CC^{n-z}, \quad \mu|_{\CC^{n-z}\times\CC^{n-z}}=\lambda, \\ \\
 \mu(X_1,X_2)=(n-2)^{\unm}(X_1+X_2)\quad ({\rm if}\;z=2),
\end{array}
$$
is a critical point of $F_n:\PP V_n\mapsto\RR$ of type $(0;n)$.  Consequently $F_n([\mu])=\frac{4}{n}$, and so
$[\mu]$ is a minimum of $F_n:\PP V_n\mapsto\RR$ by Lemma \ref{mini}.

\vs

We now study the maxima of $F_n:L_n\mapsto\RR$.  Let $\mu_{\rm he},\mu_{\rm hy}$ denote the Lie brackets defined
by
$$
\mu_{\rm he}(X_1,X_2)=X_3,\qquad \mu_{\rm hy}(X_1,X_i)=X_i, \qquad i=2,...,n,
$$
and zero otherwise.  Note that $\mu_{\rm he}$ is isomorphic to the direct sum of the 3-dimensional Heisenberg
Lie algebra and an abelian Lie algebra.

\begin{theorem}\label{max}
The functional $F_n:L_n\mapsto\RR$ attains its maximum value at $[\mu]\in L_n$ if and only if
$\mu\in\Gl(n).[\mu_{\rm he}]$.  Furthermore, $F_n([\mu_{\rm he}])=12$.
\end{theorem}

\begin{proof}
Assume that $[\mu]$ is a maximum of $F_n:L_n\mapsto\RR$.  This implies that $[\mu]$ is a critical point of
$F_n:\Gl(n).[\mu]\mapsto\RR$ and thus $[\mu]$ is a critical point of $F_n:\PP V_n\mapsto\RR$ (see Proposition
\ref{calculo}).  But then it follows from Theorem \ref{min}, (i), that $[\mu]$ is also a minimum for
$F_n:\Gl(n).[\mu]\mapsto\RR$, and hence we obtain that $F_n:\Gl(n).[\mu]\mapsto\RR$ is a constant function.
Therefore every $[\lambda]\in\Gl(n).[\mu]$ is a critical point of $F_n:V_n\mapsto\RR$ by Proposition
\ref{calculo} and so $[\mu]$ must satisfy the following rather strong condition (see Theorem \ref{min}, (ii)):
\begin{equation}\label{max1}
\Gl(n).[\mu]=\U(n).[\mu].
\end{equation}
In particular, the only possible degeneration of $\mu$ is $\mu\mapsto 0$.  By \cite[Thm 5.2]{inter} (see also
\cite{G2}), we have that $\mu_{\rm he}$ and $\mu_{\rm hy}$ are the only Lie algebras which satisfy (\ref{max1}).
It is easy to see that $[\mu_{\rm he}]$ and $[\mu_{\rm hy}]$ are critical points of type $(2<3<4;2,n-3,1)$ and
$(0<1;1,n-1)$ respectively, and so it follows from Proposition \ref{values} that
$$
F_n([\mu_{\rm he}])=12>4=F_n([\mu_{\rm hy}]),
$$
concluding the proof.
\end{proof}

We now prove the main result of this paper, that is, a description of the Lie algebras which are critical points
of $F_n$, in terms of the nilpotent critical points of the functionals $F_m$ with $m\leq n$.

\begin{theorem}\label{car}
Let $[\mu]\in L_n$ be a critical point of $F_n$ of type $\eigencero$, and consider
$\CC^n=\ggo_1\oplus\ggo_2\oplus...\oplus\ggo_r$, the eigenspace decomposition of $\Rin_{\mu}=c_{\mu}I+D_{\mu}$.
Then the following conditions hold:
\begin{itemize}
\item[(i)] $\ngo=\ggo_2\oplus...\oplus\ggo_r$ is the nilradical of $\mu$  and $\mu_{\ngo}=\mu|_{\ngo\times\ngo}$
is in his turn a critical point of the functional $F_m$ of type $(k_2<...<k_r;d_2,...,d_r)$, where $m=n-d_1$.

\item[(ii)] $\ggo_1$ is a reductive Lie subalgebra of $\mu$.

\item[(iii)] $(\ad_{\mu}{A})^*\in\Der(\mu)$ for any $A\in\ggo_1$.
\end{itemize}
Conversely, let $[\lambda]\in L_m$ be a critical point of $F_m$ of type $(k_2<...<k_r;d_2,...,d_r)$ which is
nilpotent, and let $\rg\subset\Der(\lambda)$ be a reductive Lie subalgebra of dimension $d_1$ such that
$A^*\in\Der(\lambda)$ for every $A\in\rg$. Then the semidirect product $\mu=[\cdot,\cdot]\ltimes\lambda$
determines a critical point $[\mu]\in L_n$ of $F_n$ of type $\eigencero$, where $n=d_1+m$ and $[\cdot,\cdot]$
denotes the Lie bracket of $\rg$.  The fixed inner product considered on $\CC^n=\rg\oplus\CC^m$ is the extension
of the given one on $\CC^m$ by setting
$$
\la A,B\ra=-\frac{4}{c_{\lambda}}\left(\unm\tr{\ad{A}(\ad{B})^*}+\tr{AB^*}\right),  \qquad A,B\in\rg.
$$
\end{theorem}

\begin{proof}
We first prove condition (iii), since it will be crucial in the proof of the other assertions in the theorem.
If $X\in\ggo_1$ then $[D_{\mu},\ad_{\mu}{X}]=\ad_{\mu}{D_{\mu}X}=0$, and therefore
$$
0=\tr{[D_{\mu},\ad_{\mu}{X}](\ad_{\mu}{X})^*}=\tr{D_{\mu}[\ad_{\mu}{X},(\ad_{\mu}{X})^*]}=\tr{\Rin_{\mu}[\ad_{\mu}{X},(\ad_{\mu}{X})^*]}.
$$
Thus (iii) follows from Lemma \ref{dR} (iii).

We first prove (ii).  Recall that $\ggo_1$ is a subalgebra of $\mu$ due to the $\ZZ_{\geq 0}$-gradation of $\mu$
defined by $D_{\mu}$.  Consider the orthogonal decomposition $\ggo_1=\hg\oplus\ag$, where
$\hg=\mu(\ggo_1,\ggo_1)$.  We denote by $[\cdot,\cdot]$ the Lie bracket $\mu$ restricted to the subalgebra
$\ggo_1$ and by $\ad$ its adjoint representation.  If $X\in\ggo_1$ then by (iii), the maps
$\ad{X},(\ad{X})^*:\ggo_1\mapsto\ggo_1$ have to leave invariant $\hg$, which shows that $\ag$ is an abelian
factor of $\ggo_1$ and so $[\hg,\hg]=\hg$.  Since $\ad{X}$ and $(\ad{X})^*$ must also leave invariant the center
of $\hg$ for every $X\in\ggo_1$, we obtain that the center of $\hg$ equals $0$.  Moreover, $\hg$ is semisimple.
Indeed, since $\Ker{B}=\{ H\in\hg:B(H,\cdot)\equiv 0\}$ ($B$ the Killing form of $\hg$) is
$\Der(\hg)$-invariant, we have that $(\Ker{B})^{\perp}$ is also an ideal of $\hg$.  This implies that $\Ker{B}$
is a solvable Lie algebra since its Killing form is identically zero.  But we have that
$[\Ker{B},\Ker{B}]=\Ker{B}$, and so $\Ker{B}=0$, which concludes the proof of (ii).

We now prove (i).  It follows from (ii) that $\sg=\ag\oplus\ngo$ is the radical of $\mu$.  If $A\in\ag$ belongs
to the maximal nilpotent ideal of $\sg$, then $\ad_{\mu}{A}:\ngo\mapsto\ngo$ is nilpotent.  But condition (iii)
implies that $(\ad_{\mu}{A})^*$ is a derivation, that is also nilpotent.  It is easy to see that this is
possible only if $\ad_{\mu}{A}=0$, and this implies that $A=0$.  Indeed, if $\{ X_i\}$ is an orthonormal basis
of $\ngo$ then by (\ref{defR}) we have that
$$
\begin{array}{rl}
c_{\mu}\la A,A\ra =&\la\Rin_{\mu}A,A\ra \\
=&-4\displaystyle{\sum_{ij}}\la\mu(A,X_i),X_j\ra
\overline{\la\mu(A,X_i),X_j\ra} \\
&+2\displaystyle{\sum_{ij}}\overline{\la\mu(X_i,X_j),A\ra}
\la\mu(X_i,X_j),A\ra \\
&=-4\tr{\ad_{\mu}{A}(\ad_{\mu}{A})^*}=0
\end{array}
$$
Thus $\ngo$ is the nilradical of $\mu$.

If $H\in\hg$, then by (iii) and the semisimplicity of $\hg$ we have that
$$
\ad{H}=\ad_{\mu}{H}|_{\hg}=\unm
(\ad_{\mu}{H}-(\ad_{\mu}H)^*)|_{\hg}+\unm(\ad_{\mu}{H}+(\ad_{\mu}{H})^*)|_{\hg}=\ad{X}+\ad{Y}
$$
for some $X,Y\in\hg$, and so $H=X+Y$ since $\hg$ has no center.  This implies that there is an orthonormal basis
$\{ A_i\}$ of $\ggo_1$ such that $\ad_{\mu}{A_i}|_{\ggo_1}$ is skew-hermitian for any $i$.  By (\ref{defR}) we
have that for $X,Y\in\ngo$,
$$
\begin{array}{rl}
\la\Rin_{\mu}X,Y\ra= &\la\Rin_{\mu_{\ngo}}X,Y\ra
-4\displaystyle{\sum_{ij}}\la\mu(A_i,X),X_j\ra\overline{\la\mu(A_i,Y),X_j\ra} \\
&+4\displaystyle{\sum_{ij}}\overline{\la\mu(A_i,X_j),X\ra}\la\mu(A_i,X_j),Y\ra,
\end{array}
$$
or, in other terms,
$$
\Rin_{\mu}|_{\ngo}=\Rin_{\mu_{\ngo}}+4\sum_{i}[\ad_{\mu}{A_i}, (\ad_{\mu}{A_i})^*]|_{\ngo}.
$$
Now, by applying Lemma \ref{dR},(ii) for $\mu_{\ngo}$ and using (iii), we obtain that
$$
0=\tr{\Rin_{\mu_{\ngo}}[\ad_{\mu}{A_i},(\ad_{\mu}{A_i})^*]|_{\ngo}}
$$
for all $i$.  In the same way,
$$
\begin{array}{rl}
0=&\tr{\Rin_{\mu}[\ad_{\mu}{A_i},(\ad_{\mu}{A_i})^*]}
=\tr{\Rin_{\mu}[\ad_{\mu}{A_i},(\ad_{\mu}{A_i})^*]|_{\ngo}} \\ \\
=&\tr{\Rin_{\mu}|_{\ngo}[\ad_{\mu}{A_i},(\ad_{\mu}{A_i})^*]}.
\end{array}
$$
The last two equalities imply that the hermitian operator
$$
T=\sum_i[\ad_{\mu}{A_i},(\ad_{\mu}{A_i})^*]|_{\ngo}
$$
satisfies $\tr{T^2}=0$, and so $T=0$.  Therefore $\Rin_{\mu}|_{\ngo}=\Rin_{\mu_{\ngo}}$, and consequently,
$\Rin_{\mu_{\ngo}}=c_{\mu}I+D_{\mu}|_{\ngo}$, that is, $[\mu_{\ngo}]$ is a critical point of $F_m$ of type
$(k_2<...<k_r,d_2,...,d_r)$, as asserted. This concludes the proof of (i).

\vs

We now prove the converse assertion.  Consider the subalgebra of $\Der(\lambda)$ given by
$$
\tilde{\rg}=\{ A\in\Der(\lambda):A^*\in\Der(\lambda)\}.
$$
Thus $\tilde{\rg}=\tilde{\kg}+\im\tilde{\kg}$, where
$$
\tilde{\kg}=\{ A\in\Der(\lambda):A^*=-A\},
$$
and so $\tilde{\rg}$ is reductive.  In fact, for the inner product on $\tilde{\rg}$ defined by $\la
A,B\ra_1=\tr{AB^*}$, we have that
\begin{equation}\label{car1}
(\ad{A})^{*_1}=\ad(A^*), \qquad \forall A\in\tilde{\rg},
\end{equation}
and hence $\ip_1$ is $\ad{\tilde{\kg}}$-invariant.

The Lie algebra $\tilde{\hg}=[\tilde{\rg},\tilde{\rg}]$ is semisimple and it is easy to see that each of their
simple factors $\tilde{\hg}_i$ satisfies that $A^*\in\tilde{\hg}_i$ for all $A\in\tilde{\hg}_i$.  Since
$\rg\subset\tilde{\rg}$ and so $\hg=[\rg,\rg]\subset[\tilde{\rg},\tilde{\rg}]=\tilde{\hg}$, we obtain that
$A^*\in\hg$ for any $A\in\hg$ as well, and therefore
\begin{equation}\label{car2}
\hg=\kg+\im\kg, \quad {\rm where} \quad \kg=\{ A\in\hg:A^*=-A\}.
\end{equation}
We define the fixed inner product on $\CC^n=\rg\oplus\CC^m$, also denoted by $\ip$, as follows:
$\la\rg,\CC^m\ra=0$, $\ip|_{\CC^m\times\CC^m}=\ip$ and
\begin{equation}\label{definn}
\la A,B\ra=-\frac{4}{c_{\lambda}}(\unm\tr{\ad{A}(\ad{B})^{*_1}}+\tr{AB^*}),  \qquad A,B\in\rg.
\end{equation}
It is easy to see that the adjoints of $\ad{B}$ with respect to the inner products $\ip$ and $\ip_1$ coincide.
It follows from (\ref{car1}) and (\ref{car2}) that there exists an orthonormal basis $\{ H_i\}$ of $\hg$ such
that $\ad_{\mu}{H_i}:\CC^n\mapsto\CC^n$ is skew hermitian for all $i$.  In particular, there is an orthonormal
basis $\{ A_i\}$ of $\rg$ such that $\ad_{\mu}{A_i}|_{\rg}$ is skew hermitian for all $i$.

Now, we just have to follow the proof of the first part of the theorem in the converse direction.  As in the
second part of the proof of (i), we obtain that $\Rin_{\mu}|_{\ngo}=\Rin_{\lambda}=c_{\lambda}I+D_{\lambda}$.
If $\{ X_i\}$ is an orthonormal basis of $\CC^m$ then for $A\in\rg$ and $X\in\CC^m$ we have that
$$
\la\Rin_{\mu}X,A\ra= -4\sum_{ij}\la\mu(X,X_i),X_j\ra\overline{\la\mu(A,X_i),X_j\ra}=
-4\tr{\ad_{\lambda}{X}A^*}=0,
$$
since $A^*$ is the sum of a skew-hermitian and a hermitian derivation of $\lambda$ and $\ad_{\lambda}{X}$ is
nilpotent.  If $A,B\in\rg$ and we now denote by $\{ X_i\}$ an orthonormal basis of $\CC^n$ which contains $\{
H_i\}$ , then by (\ref{defR}) we have that
\begin{equation}\label{car3}
\begin{array}{rl}
\la\Rin_{\mu}A,B\ra =& -4\displaystyle{\sum_{ij}}\la\mu(A,X_i),X_j\ra
\overline{\la\mu(B,X_i),X_j\ra} \\
&+2\displaystyle{\sum_{ij}}\overline{\la\mu(X_i,X_j),A\ra}
\la\mu(X_i,X_j),B\ra \\
=&-4\tr{\ad_{\mu}{A}(\ad_{\mu}{B})^*}+ 2\displaystyle{\sum_{ij}}\overline{\la\mu(H_i,H_j),A\ra}
\la\mu(H_i,H_j),B\ra \\
=&-2\tr{\ad{A}(\ad{B})^{*_1}}-4\tr{AB^*})
\end{array}
\end{equation}
It then follows easily from (\ref{car2}), (\ref{car3}) and the definition of $\ip|_{\rg\times\rg}$ given in
(\ref{definn}) that  $\Rin_{\mu}|_{\rg}= c_{\lambda}I$.  We conclude that
$\Rin_{\mu}=c_{\lambda}I+\tilde{D_{\lambda}}$, where $\tilde{D_{\lambda}}|_{\rg}\equiv 0$ and
$\tilde{D_{\lambda}}|_{\CC^m}=D_{\lambda}$, and hence $\mu$ is a critical point of type $\eigencero$, as it was
to be shown.
\end{proof}

From the proof of the theorem, we can deduce the following compatibility properties at a critical point $[\mu]$
of $F_n$, between the Lie bracket $\mu$ and the fixed inner product on $\CC^n$.

\begin{proposition}\label{comp}
Let $[\mu]\in L_n$ be a critical point of $F_n$ of type $\eigencero$, and consider
$\CC^n=\ggo_1\oplus\ggo_2...\oplus\ggo_r$, the eigenspace decomposition of $\Rin_{\mu}=c_{\mu}I+D_{\mu}$.  Let
$\ggo_1=\hg\oplus\ag$ with $\hg$ semisimple and $\ag$ the center of $\ggo_1$.  Then the following conditions
hold:
\begin{itemize}
\item[(i)] $\ad_{\mu}{A}$ is a normal operator for every $A\in\ag$.

\item[(ii)] The real subalgebra $\kg=\{ A\in\hg:(\ad_{\mu}A)^*=-\ad_{\mu}A\}$ is a maximal compact subalgebra of
$\hg$, that is, $\hg=\kg+\im\kg$.

\item[(iii)] The hermitian inner product on $\ggo_1$ satisfies
$$
\la A,B\ra=-\frac{4}{c_{\mu}}\left(\unm\tr{\ad_{\mu}{A}
(\ad_{\mu}{B})^*}|_{\hg}+\tr{\ad_{\mu}{A}(\ad_{\mu}{B})^*}|_{\ngo}\right), \qquad A,B\in\ggo_1.
$$
\end{itemize}
\end{proposition}

\begin{proof}
Parts (ii) and (iii) follow directly from the proof of Theorem \ref{car}. We now prove (i).  If $A\in\ag$ then
$(\ad_{\mu}{A})^*\in\Der(\mu)$ by Theorem \ref{car} (iii).  This implies that
$$
0=\ad_{\mu}((\ad_{\mu}{A})^*A)=[(\ad_{\mu}{A})^*,\ad_{\mu}{A}],
$$
which concludes the proof.
\end{proof}

Roughly speaking, Theorem \ref{car} says that the study of the Lie algebras which are critical points of $F_n$
reduces to the understanding of those which are nilpotent.  In other words, it suffices to describe the
Kirwan-Ness quotient $\nca_n\kir\Gl(n)$, where $\nca_n\subset\lca_n\subset\PP V_n$ is the algebraic subvariety
of nilpotent Lie algebras.

\begin{remark}\label{riem}
{\rm We would like to point out a somewhat mysterious characterization of the Kirwan-Ness quotient
$\nca_n\kir\Gl(n)$ in terms of Riemannian geometry: the classification of the nilpotent critical points of
$F_n:L_n\mapsto\RR$ is equivalent to the classification up to isometry of all the left invariant Riemannian
metrics on nilpotent Lie groups which are Ricci solitons (see \cite{soliton}), and also of all the Einstein left
invariant Riemannian metrics on solvable Lie groups (see \cite{H, critical}). A crucial point here is that the
moment map $\Rin_{\mu}$ coincides with the Ricci operator of a certain Riemannian metric naturally associated
with $\mu$.}
\end{remark}

We now give some known examples of nilpotent critical points, most of which come from the interplay mentioned
above.

\begin{example}\label{exanil}
{\rm (i) A two-step nilpotent Lie algebra $\mu$ is a critical point of $F_n$ of type $(1<2;d_1,d_2)$ if and only
if $\Rin_{\mu}|_{\zg(\mu)}$ and $\Rin_{\mu}|_{\zg(\mu)^{\perp}}$ are both a multiple of the identity, where
$\zg(\mu)$ is the center of $\mu$ and $\dim{\zg(\mu)}=d_2$.  Thus any Heisenberg-type Lie algebra (see
\cite{BTV}) and any two-step nilpotent Lie algebra constructed via a representation of a compact Lie group (see
\cite{rep}), contain a critical point in their orbits (over $\RR$).  In \cite{GT}, curves of critical points of
type $(1<2;5,5)$ and $(1<2;6,3)$ are given.

\no (ii) Any nilpotent Lie algebra with a codimension one abelian ideal contains a critical point of $F_n$ in
its orbit (see $\S$\ref{abelian}).

\no (iii) Every nilpotent Lie algebra of dimension $\leq 5$ is isomorphic to a critical point of $F_n$ (see
\cite{L1}).

\no (iv) The lowest possible dimension for the existence of a curve of non-isomorphic nilpotent critical points
of $F_n$ is $n=7$, and an example of such a curve is given in \cite{L1}.

\no (v) It has been recently proved in \cite{Wll} that any $6$-dimensional nilpotent Lie algebra contains a
critical point of $F_6$ in its orbit. }
\end{example}

On the other hand, since any nilpotent critical point is necessarily $\NN$-graded, we have that a
characteristically nilpotent Lie algebra (i.e. $\Der(\mu)$ nilpotent) can never be a critical point of $F_n$.
The first dimension where these Lie algebras appear is $n=7$.

\section{Lie algebras with a codimension one abelian ideal}\label{abelian}

We consider the closed subset of $\lca_{n+1}$ given by
$$
\aca=\{\mu\in\lca_{n+1}:\mu\;\mbox{has an abelian ideal of dimension $n$}\}.
$$
The aim of this section is to illustrate, via a study of $\aca$, most of the notions considered in this paper.
Since the dimension of the maximal abelian subalgebra increases with a degeneration, we have that the stratum
$S_{(0<1;1,n)}\subset\aca$ (see (\ref{stratum})).

By fixing a decomposition $\CC^{n+1}=\CC H\oplus\CC^n$, every element in $\aca$ is isomorphic to a $\mu\in\aca$
having $\CC^n$ as the required abelian ideal.  In this case, $\mu$ is determined by the matrix
$A=\ad_{\mu}{H}|_{\CC^n}\in\glg(n)$, and so it will be denoted by $\mu_A$.  It is easy to see that $\mu_A$ is
isomorphic to $\mu_B$, $A,B\in\glg(n)$, if and only if $A$ is conjugate to $B$ up to scaling.  In other words,
isomorphism classes of $\aca$ are parameterized by $\glg(n)/\Gl(n)$, under the action
\begin{equation}\label{abiso}
\vp.A=(\det{\vp})\vp A\vp^{-1}, \qquad \vp\in\Gl(n),\;A\in\glg(n).
\end{equation}
It is not difficult to prove that for a non-nilpotent $A$,
\begin{equation}\label{abder1}
\Der(\mu_A)=\left\{ \left[\begin{array}{c|c}
0&0\\
\hline \ast&B
\end{array}\right]:\quad BA=AB, \quad B\in\glg(n)\right\},
\end{equation}
and if $A$ is nilpotent then
\begin{equation}\label{abder2}
\Der(\mu_A)=\CC D\oplus\left\{ \left[\begin{array}{c|c}
0&0\\
\hline \ast&B
\end{array}\right]:\quad BA=AB,\quad B\in\glg(n)\right\},
\end{equation}
where
$$
D=\left[\begin{array}{c|c}
1&0\\
\hline 0&D_1
\end{array}\right],\qquad [D_1,A]=A.
$$
Using the formula (\ref{defR}) as in the proof of Theorem \ref{car}, we can  calculate the moment map and the
functional $F_{n+1}$ on $\aca$.  We assume that the decomposition $\CC^{n+1}=\CC H\oplus\CC^n$ is orthogonal and
$||H||=1$.

\begin{proposition}\label{abR}
For any $A\in\glg(n)$ we have that
$$
\Rin_{\mu_A}=\left[\begin{array}{c|c}
-4\tr{AA^*}&0\\
\hline 0&4[A,A^*]
\end{array}\right].
$$
Consequently, $||\mu_A||^2=2\tr{AA^*}$ and $F_{n+1}([\mu_A])=4+16\tr{[A,A^*]^2}$.
\end{proposition}

We deduce from the formula of $F_{n+1}([\mu_A])$ that $F_{n+1}$ measures how far is $A$ from being normal, and
so if $A$ is semisimple then the orbit $\Gl(n+1).[\mu_A]$ will contain a critical point of $F_{n+1}$ (see
Theorem \ref{min}).  Indeed, for any semisimple $A$ there exists $\vp\in\Gl(n)$ such that $\vp A\vp^{-1}$ is
normal with respect to $\ip$, and hence if $\psi\in\Gl(n+1)$ is defined by $\psi|_{\CC H}=1$ and
$\psi|_{\CC^n}=\vp$, then $\psi.\mu_A=\mu_{\vp A\vp^{-1}}$, and so $\psi.\mu_A$ is a critical point of type
$(0<1;1,n)$.  However as a counterpart, there is also a critical point of $F_{n+1}$ in the orbits of the
nilpotent $\mu_A$'s.

\begin{proposition}\label{abcrit}
The orbit $\Gl(n+1).[\mu_A]$ contains a critical point of $F_{n+1}$ if and only if $A$ is semisimple or
nilpotent. If $A$ is not nilpotent, then $[\mu_A]$ is a critical point of $F_{n+1}$ if and only if $A$ is
normal, and in this case, it is of type $(0<1;1,n)$.
\end{proposition}

\begin{proof}
We first assume that $A$ is not nilpotent.  If $[\mu_A]$ is a critical point of $F_{n+1}$ then
$\Rin_{\mu_A}\in\RR I\oplus\Der(\mu_A)$ (see Proposition \ref{calculo},(iii)).  Thus it follows from
(\ref{abder1}) and Proposition \ref{abR} that $[A,A^*]$ commutes with $A$, and hence
$$
\left[\begin{array}{c|c}
0&0\\
\hline 0&[A,A^*]
\end{array}\right]=
\left[\left[\begin{array}{c|c}
0&0\\
\hline 0&A
\end{array}\right],
\left[\begin{array}{c|c}
0&0\\
\hline 0&A^*
\end{array}\right]\right]
$$
is a derivation of $\mu_A$.  Now, parts (ii) and (iii) of Lemma \ref{dR} imply that
$$
\left[\begin{array}{c|c}
0&0\\
\hline 0&A^*
\end{array}\right]
$$
is also a derivation of $\mu_A$, concluding that $A$ is normal (see (\ref{abder1})).  The type of $[\mu_A]$ can
be deduced from the formula
$$
\Rin_{\mu_A}=-4\tr{AA^*}I+ \left[\begin{array}{c|c}
0&0\\
\hline 0&4\tr{AA^*}I
\end{array}\right].
$$
We now consider $A$ nilpotent.  For each $r$-tuple $(n_1,...,n_r)$ of integer numbers satisfying $n_1\geq...\geq
n_r\geq 0$ and $(n_1+1)+...+(n_r+1)=n$, we consider the $n\times n$-matrix $A_{(n_1,...,n_r)}$ obtained by the
direct sum of the $r$ blocks of the form
$$
\left[\begin{array}{cccccc}
0&&&&&\\
n_i^{\unm}&\ddots&&&&\\
&\ddots&0&&&\\
&&(jn_i-j(j-1))^{\unm}&0&&\\
&&&\ddots&\ddots&\\
&&&&n_i^{\unm}&0
\end{array}\right], \qquad {\rm where}\; j=1,...,n_i.
$$
By a straightforward computation, one can see that $[\mu_{A_{(n_1,...,n_r)}}]$ is a critical point of $F_{n+1}$
for any $(n_1,...,n_r)$.  If $A\in\glg(n)$ is nilpotent then $\mu_A$ is isomorphic to
$[\mu_{A_{(n_1,...,n_r)}}]$, where $n_1+1\geq...\geq n_r+1$ are the dimensions of the Jordan blocks of $A$. Thus
every orbit $\Gl(n+1).[\mu_A]$ with $A$ nilpotent contains a critical point of $F_{n+1}$.
\end{proof}

Therefore, the Kirwan-Ness quotient is given by
$$
\aca\kir\Gl(n+1)=\{\mu_A:A\;\mbox{semisimple or nilpotent}\}=\PP\CC^n/S_n\cup\bigcup\{\mu_{A_{(n_1,...,n_r)}}\}
$$
for $n_1\geq...\geq n_r\geq 0$, $n_1+...+n_r=n-r$ (see the proof of the above proposition), where $S_n$ denotes
the symmetric group permuting the coordinates of $\PP\CC^n$.

Consider for each $A\in\glg(n)$ the decomposition $A=S+N$ in their semisimple and nilpotent parts.  If $S\ne 0$,
it is easy to see that the $-\grad(F_{n+1})$-flow starting from the point $[\mu_{A}]$ converges to the critical
point $[\mu_S]$.  Thus the stratum of type $(0<1;1,n)$ is given by
$$
S_{(0<1;1,n)}=\bigcup_{S\ne 0}\Gl(n+1).\mu_A =\{\mu\in\aca:\mu\;\mbox{is not nilpotent}\}.
$$
  It is proved in
\cite{L1} that the type of the critical point  $[\mu_{A_{(n_1,...,n_r)}}]$ is given by
\begin{itemize}
\item $(1<\theta-\frac{n_1}{2}<\theta-\frac{n_1}{2}+1<...<...<\theta+\frac{n_1}{2};1,...)$, \quad if
$n_i\equiv\epsilon\;(2)$, for each $i$,

\item $(2<3<4;2,n-3,1)$, \quad if $n_1=1$, $n_i=0$ for each $i\geq 2$,

\item $(2<2\theta-n_e<2\theta-n_e+2<...<2\theta+n_o;1,...)$, \quad otherwise,
\end{itemize}
where $n_e$ and $n_o$ are the greatest even and odd numbers among the $n_i's$ respectively and
$$
\theta=1+\frac{n_1(n_1+1)(n_1+2)+...+n_r(n_r+1)(n_r+2)}{12}.
$$
It is easy to check that the type of $[\mu_{A_{(n_1,...,n_r)}}]$ coincides with the type of
$[\mu_{A_{(n_1',...,n_r')}}]$ if and only if $n_i=n_i'$ for all $i$.  This implies that the other strata $S_{
\alpha}\cap\aca$ are just the orbits $\Gl(n+1).[\mu_A]$ for each nilpotent $A$.

\section{Critical points and closed orbits}\label{critclosed}

In this section, we show how the critical points of $F_n$ of a given type can be viewed (up to isomorphism) as
the categorical quotient of a suitable action (see \cite[Section 9]{Ns} for the general case).  We fix a type
$\alpha=\eigen$, and denote by
$$
V_{\alpha}=\{\mu\in V_n:D_{\alpha}\in\Der(\mu)\}, \qquad
G_{\alpha}=Z_{\Gl(n)}(D_{\alpha})=\Gl(d_1)\times...\times\Gl(d_r),
$$
where
$$
D_{\alpha}= \left[\begin{array}{ccc}
k_1I_{d_1}&&\\
&\ddots&\\
&&k_rI_{d_r}
\end{array}\right],
$$
and $I_{d_i}$ denotes the $d_i\times d_i$ identity matrix.  We consider the reductive subgroup of $G_{\alpha}$
defined by
$$
\tilde{G}_{\alpha}=\left\{ g\in G_{\alpha}:\prod_{i=1}^{r}(\det{g_i})^{k_i}=\det{g}=1\right\}.
$$
We are interested in the action of $\tilde{G}_{\alpha}$ on $V_{\alpha}$.  The corresponding Lie algebras satisfy
$$
\ggo_{\alpha}=\tilde{\ggo}_{\alpha}\oplus\CC I\oplus\CC D_{\alpha}, \qquad \tilde{\ggo}_{\alpha}=\{
A\in\ggo_{\alpha}:\tr{AD_{\alpha}}=\tr{A}=0\}.
$$
Let $\kg_{\alpha},\tilde{\kg}_{\alpha}$ denote the Lie algebras of the maximal compact subgroups.  It is easy to
see that the moment map $m:\PP V_{\alpha}\mapsto\im\kg_{\alpha}$ corresponding to the action of $G_{\alpha}$ on
$V_{\alpha}$ is just the restriction of $m:\PP V\mapsto\im\ug(n)$ to $\PP V_{\alpha}$, and the moment map
$\tilde{m}:\PP V_{\alpha}\mapsto\im\tilde{\kg}_{\alpha}$ of the $\tilde{G}_{\alpha}$-action on $V_{\alpha}$
coincides with the composition of $m:\PP V_{\alpha}\mapsto\im\kg_{\alpha}$ with the orthogonal projection
$\im\kg_{\alpha}\mapsto\im\tilde{\kg}_{\alpha}$.

Recall that $[\mu]\in\PP V_{\alpha}$ is a critical point of type $\alpha$ if and only if
$m(\mu)=\Rin_{\mu}\in\CC(c_{\alpha}I+D_{\alpha})$, which is equivalent to $\Rin_{\mu}\in\CC I\oplus\CC
D_{\alpha}$.  Since $\im\kg_{\alpha}=\im\tilde{\kg}_{\alpha}\oplus(\CC I+\CC D_{\alpha})$ is an orthogonal
decomposition, it follows that
$$
\tilde{m}^{-1}(0)=\PP V_{\alpha}\cap C_{\alpha},
$$
where $C_{\alpha}$ is the set of critical points of $F_n$ of type $\alpha$.  This means that the orbit
$\tilde{G}_{\alpha}.\mu$ of $\mu\in V_{\alpha}$ is closed if and only if $\tilde{G}_{\alpha}.[\mu]$ intersects
$C_{\alpha}$.  The open set of semistable points is precisely
\begin{equation}\label{ssstratum}
\PP V_{\alpha}^{\rm ss}=\PP V_{\alpha}\cap S_{\alpha},
\end{equation}
and the categorical quotient $\PP V_{\alpha}//\tilde{G}_{\alpha}$ coincides with
$$
\PP V_{\alpha}\cap C_{\alpha}/\tilde{K_{\alpha}}=C_{\alpha}/\U(n),
$$
which parameterizes the set of critical points of type $\alpha$ up to isomorphism.  This shows that
$C_{\alpha}/\U(n)$ is a projective algebraic variety.  If $\alpha_1,...,\alpha_s$ denote the different types of
critical points of $F_n:L_n\mapsto\RR$ then the stratification of the Kirwan-Ness quotient described in Section
\ref{ness} is given by
$$
L_n\kir\Gl(n)=X_{\alpha_1}\cup...\cup X_{\alpha_s} \qquad \mbox{(disjoint union},
$$
where each $X_{\alpha_i}$ is homeomorphic to $C_{\alpha_i}/\U(n)$ and the following frontier property holds:
$$
\overline{X}_{\alpha}\subset X_{\alpha}\cup\bigcup_{\beta>\alpha}X_{\beta}.
$$
Recall that $\alpha<\beta$ if and only if the corresponding critical values satisfy $F_n(\alpha)<F_n(\beta)$
(see Proposition \ref{values}).

{\small
\begin{table}\label{isofour}
$$
\begin{array}{lc}
\hline \\
\ggo & \mbox{Lie brackets} \\ \\
\hline\\
\CC^4 & --- \\ \\
\ngo_3\oplus\CC & [x_1,x_2]=x_3 \\ \\
\rg_2\oplus\CC^2 & [x_1,x_2]=x_1 \\ \\
\rg_3\oplus\CC & [x_1,x_2]=x_2, [x_1,x_3]=x_2+x_3 \\ \\
\rg_{3,\lambda}\oplus\CC & [x_1,x_2]=x_2, [x_1,x_3]=\lambda x_3, \lambda\in\CC, \; 0<|\lambda|\leq 1 \\ \\
\rg_2\oplus\rg_2 & [x_1,x_2]=x_1, [x_3,x_4]=x_3 \\ \\
\slg_2\oplus\CC & [x_1,x_2]=x_3, [x_1,x_3]=-2x_1, [x_2,x_3]=2x_2 \\ \\
\ngo_4 & [x_1,x_2]=x_3, [x_1,x_3]=x_4 \\ \\
\ggo_1(\alpha) & [x_1,x_2]=x_2, [x_1,x_3]=x_3, [x_1,x_4]=\alpha x_4, \alpha\in\CC^* \\ \\
\ggo_2(\alpha,\beta) & [x_1,x_2]=x_3, [x_1,x_3]=x_4, [x_1,x_4]=\alpha x_2-\beta x_3+x_4, \\
& \alpha\in\CC^*, \beta\in\CC \; {\rm or}\; \alpha,\beta=0 \\ \\
\ggo_3(\alpha) & [x_1,x_2]=x_3, [x_1,x_3]=x_4, [x_1,x_4]=\alpha(x_2+x_3), \alpha\in\CC^* \\ \\
\ggo_4 & [x_1,x_2]=x_3, [x_1,x_3]=x_4, [x_1,x_4]=x_2 \\ \\
\ggo_5 & [x_1,x_2]=\unt x_2+x_3, [x_1,x_3]=\unt x_3, [x_1,x_4]=\unt x_4 \\ \\
\ggo_6 & [x_1,x_2]=x_2, [x_1,x_3]=x_3, [x_1,x_4]=2x_4, [x_2,x_3]=x_4 \\ \\
\ggo_7 & [x_1,x_2]=x_3, [x_1,x_3]=x_2, [x_2,x_3]=x_4 \\ \\
\ggo_8(\alpha) & [x_1,x_2]=x_3, [x_1,x_3]=-\alpha x_2+x_3, [x_1,x_4]=x_4, [x_2,x_3]=x_4, \alpha\in\CC \\
\\ \hline
\end{array}
$$
\caption{Complex Lie algebras of dimension $4$}
\end{table} }

\begin{example}
{\rm Consider the type $\alpha=(0;n)$.  We have that $V_{\alpha}=V_n$, $G_{\alpha}=\Gl(n)$ and
$\tilde{G}_{\alpha}=\Sl(n)$. Thus the set of semistable points $\PP V_n^{\rm ss}=S_{\alpha}$ is open in $\PP
V_n$ and consequently $S_{\alpha}\cap L_n$ is open in $L_n$.  We know that  $S_{\alpha}\cap L_n$ is precisely
the set of $n$-dimensional semisimple Lie algebras (see Remark \ref{sson}).  If $S_{\alpha}\cap
L_n=\Gl(n).[\mu_1]\cup...\cup\Gl(n).[\mu_m]$, where $\mu_1,...,\mu_m$ are the $n$-dimensional semisimple Lie
algebras up to isomorphism, then it follows from the fact that $\dim{\Der(\mu_j}=n$ for all $j$ that
$\overline{\Gl(n).[\mu_j]}\cap S_{\alpha}\cap  L_n=\Gl(n).[\mu_j]$, and so $\Gl(n).[\mu_j]$ is closed in
$S_{\alpha}\cap L_n$ for any $j$.  This implies that $\Gl(n).[\mu_i]$ is open in $S_{\alpha}\cap L_n$ and so in
$L_n$ for any $i$.  As an immediate consequence, $\Gl(n).\mu_i$ is open in $\lca_n$ for any $i$.  We therefore
obtain an alternative proof of the rigidity of semisimple Lie algebras, which makes use of any cohomology
argument.

On the other hand, we also conclude from Remark \ref{sson} that an orbit $\Sl(n).\mu\subset\lca_n$ is closed if
and only if $\mu$ is semisimple, and  any non-semisimple Lie algebra $\lambda$ is in the null-cone.  Thus the
categorical quotient is $L_n//\Sl(n)=\{\mu_1,...,\mu_m\}$. }
\end{example}

We now summarize the results obtained in the above example, as a corollary of Theorem \ref{ss}.

\begin{corollary} {\rm (i)} The $\Gl(n)$-orbit of any semisimple Lie algebra is open in $\lca_n$.

{\rm (ii)} The $\Sl(n)$-orbit of $\mu\in\lca_n$ is closed if and only if $\mu$ is semisimple.

{\rm (iii)} $0\in\overline{\Sl(n).\lambda}$ for any non-semisimple $\lambda\in\lca_n$.
\end{corollary}

\begin{example}
{\rm If $\alpha=(0<1;1,n)$, then $V_{\alpha}=\{\mu_A:A\in\glg(n)\}\subset\aca\subset\lca_{n+1}$ (see
$\S$\ref{abelian}), $G_{\alpha}=\CC^*\times\Gl(n)$ and $\tilde{G}_{\alpha}=\{ 1\}\times\Sl(n)$.  The semistable
points are
$$
V_{\alpha}^{\rm ss}=V_{\alpha}\cap S_{\alpha}=\{\mu_A:A\;\mbox{is not nilpotent}\}
$$
and the null-cone $N=\{\mu_A:A\;\mbox{is nilpotent}\}$.  An orbit $\tilde{G}_{\alpha}.\mu_A$ is closed if and
only if $A$ is semisimple, and the categorical quotient is given by
$$
\PP V_{\alpha}//\tilde{G}_{\alpha}=\PP \CC^n/S_n,
$$
which parameterizes diagonal matrices up to conjugation and scaling. }
\end{example}

In what follows, we completely develop the case $n=4$.  A list of all $4$-dimensional complex Lie algebras up to
isomorphism is given in Table 1 (see \cite{BS}).  In Table 2 we give the stratification for $L_4$.  For each
type $\alpha$, we describe the stratum $S_{\alpha}$ and the categorical quotient
$L_4//\tilde{G}_{\alpha}=C_{\alpha}/\U(4)$, as well as the critical values.  We also denote by $\ggo$ the
$\Gl(n)$-orbit in $L_4$ of Lie algebras isomorphic to $\ggo$.

\begin{table}\label{strfour}
$$
\begin{array}{c|c|c|c}
\hline
\alpha\;\mbox{(type)} & S_{\alpha} & C_{\alpha}/\U(4) & F_4(C_{\alpha}) \\

\hline &&&\\

(0<1;3,1)& \slg_2\oplus\CC & \{ p\} & \frac{4}{3} \\

\hline  &&\\

(0<1;2,2) & \rg_2\oplus\rg_2 & \{ p\} & 2 \\

\hline  &&\\

(0<1<2;1,2,1) & \mbox{solvable with nilradical}\;\ngo_3 & \CC\cup\{ p\} & 3\\

\hline  &&\\

(0<1;1,3) & \mbox{solvable with nilradical}\; \CC^3 & \PP \CC^3/S^3 & 4\\

\hline  &&\\

(1<2<3<4;1,1,1,1)& \ngo_4 & \{ p\} & 6\\

\hline  &&\\

(2<3<4;2,1,1) & \ngo_3\oplus\CC & \{ p\} & 12\\

\hline
\end{array}
$$
\caption{Stratification for $L_4$}
\end{table}

For $\alpha=(0<1<2;1,2,1)$, we have that $S_{\alpha}=\{ \ggo_8(c):c\in\CC\}\cup\{\ggo_6,\ggo_7\}$.  All of these
orbits meet $C_{\alpha}$ excepting $\ggo_8(\frac{1}{4})$, for which the negative gradient flow of $F_4$
converges to the critical point in $\ggo_6$. Recall that the case $\alpha=(0<1;1,3)$ has been studied in
$\S$\ref{abelian}.  Thus the Kirwan-Ness quotient $L_4\kir\Gl(4)$ classifies all $4$-dimensional Lie algebras
except for $\ggo_8(\frac{1}{4})$, $\ggo_3(\frac{27}{4})$, $\ggo_5$, $\ggo_2(\frac{1}{27},\frac{1}{3})$ and
$\ggo_2(\frac{\gamma}{(\gamma+2)^3}, \frac{2\gamma+1}{(\gamma+2)^2})$, $\gamma\in\CC\setminus\{ 2\}$.  The
negative gradient flow of $F_4$ starting from $\ggo_3(\frac{27}{4})$, $\ggo_5$,
$\ggo_2(\frac{1}{27},\frac{1}{3})$ and $\ggo_2(\frac{\gamma}{(\gamma+2)^3}, \frac{2\gamma+1}{(\gamma+2)^2})$
converges to the critical point in $\ggo_1(-2)$, $\ggo_1(1)$, $\ggo_1(1)$ and $\ggo_1(\gamma)$, respectively.

The classification of all possible degenerations for $4$-dimensional complex Lie algebras given in \cite{BS} can
be used to check the validity of the frontier property of the stratifications of $L_4$ and $L_4\kir\Gl(n)$.

We note that the closure $\overline{S_{\alpha}}$ of any non-nilpotent stratum $S_{\alpha}$ (i.e. $k_1=0$) gives
rise an irreducible component of $\lca_4$.  One can see that this is true for every $\lca_n$ with $n\leq 7$, by
using the results given in \cite{CD}.  We do not know if the irreducible components of $\lca_n$ can be described
in this way for any value of $n$.

\begin{example}
{\rm The approach considered in this section has been used successfully by L. Galitski and D. Timashev \cite{GT}
in the study of two-step nilpotent Lie algebras.  For a type of the form $\alpha=(1<2;d_1,d_2)$ one gets
$V_{\alpha}=\Lambda^2(\CC^{d_1})^*\otimes\CC^{d_2}$ and $\tilde{G}_{\alpha}=\Sl(d_1)\times\Sl(d_2)$.  Thus the
closed $\tilde{G}_{\alpha}$-orbits in $V_{\alpha}$ are precisely those which contain a critical point of type
$\alpha$.  In \cite{GT} the quotient space $V_{\alpha}/\tilde{G}_{\alpha}$ is studied by using methods in
geometric invariant theory, in the cases $(d_1,d_2)=(5,5)$ and $(6,3)$.  This allows to complete the
classification of two-step nilpotent Lie algebras of dimension $\leq 9$.  They exploit the fact that in these
cases, the action of $\tilde{G}_{\alpha}$ on $V_{\alpha}$ is not only visible (i.e. $p^{-1}(p(\mu))$ contains
finitely many orbits for any $\mu$, where $p:V_{\alpha}\mapsto V_{\alpha}//\tilde{G}_{\alpha}$), but it is also
a $\theta$-group.  }
\end{example}

\begin{remark}
{\rm As a consequence of Theorem \ref{car}, we have that any nilpotent critical point of $F_n$ admits an
$\NN$-gradation.  We do not know if the converse assertion is true.  For instance, if the representations
$\Lambda^2\Gl(d_1)^*\otimes\Gl(d_2)$ with $d_1+d_2=n$ were nice enough, in the sense that any orbit contains a
critical point of $F_n$, then we would be able to describe the moduli space of all $n$-dimensional two-step
nilpotent Lie algebras up to isomorphism as a finite union of categorical quotients.  However, we have reasons
to believe that this would be too optimistic for large values of $n$.  }
\end{remark}

\end{document}